\documentclass{aero}

\usepackage[utf8]{inputenc}
\usepackage{textcomp}
\usepackage[version=4]{mhchem}
\usepackage{siunitx}
\usepackage{lipsum}

\usepackage[nonumberlist,acronym]{glossaries} 
\usepackage{graphicx,setspace,parskip,bm,amsmath,mathrsfs} 
\usepackage{longtable,tabularx,multicol,multirow,threeparttablex} 
\usepackage{float} 
\usepackage{makecell} 
\usepackage{stackengine}
\usepackage{booktabs} 
\usepackage{upgreek}
\usepackage{footnpag,marginnote} 
\usepackage{booktabs} 
\usepackage[font=small]{caption} 
\usepackage{subcaption} 
\usepackage{comment}
\usepackage{microtype} 
\usepackage{makecell} 
\usepackage{setspace} 
\usepackage{pifont}
\usepackage{xcolor} 
\usepackage{colortbl} 
\usepackage{mathtools}
\usepackage{centernot}
\usepackage{pgfplots,bm}
\pgfplotsset{compat=newest}
\pgfplotsset{plot coordinates/math parser=false}
\newlength\figureheight
\newlength\figurewidth
\usepackage{environ,tikz} 
\usetikzlibrary{mindmap, trees, calc}
\usetikzlibrary{positioning,fit}
\usetikzlibrary{through}
\usetikzlibrary{shapes.geometric, arrows}
\usetikzlibrary{automata,positioning}
\usepgfplotslibrary{groupplots}
\usepackage{natbib}

\tikzset{block/.style={draw,thick,text width=2cm,minimum height=1cm,align=center},
         line/.style={-latex}
}
\usepackage{etoolbox}
\usepackage{algorithm}
\usepackage{algpseudocode}

\definecolor{phdgray}{HTML}{7f7f7f}
\definecolor{phdred}{HTML}{E57373}

\definecolor{phdred2}{HTML}{E63946}
\definecolor{phdblue}{HTML}{4a90e2}
\definecolor{phdgreen}{HTML}{4caf50}
\definecolor{phdpurple}{HTML}{6B4F9E}
\definecolor{phdgold}{HTML}{f5a623}
\definecolor{phdcyan}{HTML}{00bfae}

\DeclareSIUnit{\AU}{AU}
\DeclareSIUnit{\day}{d}
\DeclareSIUnit{\deg}{deg}
\DeclareSIUnit{\year}{yr}

\newcommand{\Eq}[1]{Eq.~\eqref{#1}}

\newcommand{\Fig}[1]{Fig.~\ref{#1}}

\newcommand{\Tab}[1]{Table~\ref{#1}}
\newcommand{\Alg}[1]{Alg.~\ref{#1}}
\newcommand{\Sec}[1]{Section~\ref{#1}}

\newcommand{\tr}{\textcolor{black}}
\newcommand{\tb}{\textcolor{black}}
\newcommand{\tp}{\textcolor{black}}
\newcommand{\tg}{\textcolor{black}}

\newcommand{\Mnota}[2]{\marginpar{\footnotesize{\tr{#1} \\ \tb{#2}}}} 

\newcommand{\ie}{i.\,e.,~}

\newcommand{\lastdate}{Oct 15, 2025}

\newacronym{ABM}{ABM}{Adams-Bashforth-Moulton}
\newacronym{ADS}{ADS}{automatic domain splitting}
\newacronym{AUL}{AUL}{augmented Lagrangian}
\newacronym{CR3BP}{CR3BP}{circular restricted three-body problem}
\newacronym{EMS}{EMS}{Earth—Moon system}
\newacronym{STM}{STM}{state-transition matrix}
\newacronym{TPBVP}{TPBVP}{two-point boundary value problem}
\newacronym{DRO}{DRO}{distant retrograde orbit}
\newacronym{NRHO}{NRHO}{near-rectilinear halo orbit}
\newacronym{DROs}{DROs}{distant retrograde orbits}
\newacronym{NRHOs}{NRHOs}{near-rectilinear halo orbits}
\newacronym{DA}{DA}{differential algebra}
\newacronym{DACE}{DACE}{differential algebra core engine}
\newacronym{DADDy}{DADDy}{differential algebra-based differential dynamic programming}
\newacronym{DDP}{DDP}{differential dynamic programming}
\newacronym{DNC}{DNC}{did not converge}
\newacronym{DOPRI8}{DOPRI8}{Dormand-Prince 8th-order embedded Runge-Kutta method}
\newacronym{DOP853}{DOP853}{Dormand-Prince of order 8(5,3) embedded Runge-Kutta method}
\newacronym{EOM}{EoM}{equations of motion}
\newacronym{GEO}{GEO}{geostationary orbit}
\newacronym{GTO}{GTO}{geostationary transfer orbit}
\newacronym{HDDP}{HDDP}{hybrid differential dynamic programming}
\newacronym{IC}{IC}{initial conditions}
\newacronym{iLQR}{iLQR}{iterative linear-quadratic regulator}
\newacronym{JPL}{JPL}{Jet Propulsion Laboratory}
\newacronym{MB}{MB}{megabytes}
\newacronym{LEO}{LEO}{low Earth orbit}
\newacronym{MEO}{MEO}{medium Earth orbit}
\newacronym{NASA}{NASA}{National Aeronautics and Space Administration}
\newacronym{NEA}{NEA}{near-Earth asteroid}
\newacronym{NEO}{NEO}{near-Earth object}
\newacronym{NSG}{NSG}{non-spherical gravity}
\newacronym{ODE}{ODE}{ordinary differential equation}
\newacronym{PECE}{PECE}{predictor-corrector}
\newacronym{PCK}{PCK}{planetary constants kernel}
\newacronym{PDS}{PDS}{Planetary Data System}
\newacronym{PN}{PN}{projected Newton}
\newacronym{RAAN}{RAAN}{right ascension of the ascending node}
\newacronym{RHS}{RHS}{right-hand side}
\newacronym{RK}{RK}{Runge-Kutta}
\newacronym{RK78}{RK78}{order 7(8) Runge–Kutta–Fehlberg embedded method}
\newacronym{RSS}{RSS}{root sum square}
\newacronym{RTN}{RTN}{radial-tangential-normal}
\newacronym{RPF}{RPF}{roto-pulsating frame}
\newacronym{SDDP}{SDDP}{stochastic differential dynamic programming}
\newacronym{STTs}{STTs}{state transition tensors}
\newacronym{ToF}{ToF}{time-of-flight}
\newacronym{UOM}{uom}{unit of measurement}
\newacronym{VSVO}{VSVO}{variable-step variable-order}
\newacronym{VV}{V\&V}{verification and validation}
\newacronym{DME}{DME}{direct mapping evaluation}
\newacronym{MSOP}{MSOP}{mapping-seeded optimisation with polynomials}
\newacronym{OP}{OP}{optimisation with polynomials}
\newacronym{PWO}{PWO}{point-wise optimisation}

\ADsetup{
title    = {Taylor polynomial-based constrained solver for fuel-optimal low-thrust trajectory optimisation},
author   = {Thomas Caleb$^{1}$, Roberto Armellin$^{2}$\cor{}, Spencer Boone$^{1}$, and St\'ephanie Lizy-Destrez$^{1}$},
address  = {
	1. DCAS/SaCLaB, ISAE-SUPAERO, Toulouse 31055, France \sep
	2. Te P\=unaha \=Atea -- Space Institute, University of Auckland, Auckland 1010, New Zealand
},
email    = {Roberto Armellin, roberto.armellin@auckland.ac.nz},
abstract = {This paper presents \gls*{DADDy}, a publicly available C++ framework for constrained, fuel-optimal low-thrust trajectory optimisation.
The method uses \gls*{DA} for two purposes: automatic differentiation and high-order Taylor expansions of the dynamics.
These expansions replace many expensive numerical propagations with polynomial evaluations, reducing computational effort while preserving solution quality.
The framework combines two complementary modules. First, a \gls*{DDP}/\gls*{iLQR} stage computes an almost-feasible trajectory from imperfect initial guesses.
Second, a polynomial-accelerated Newton stage enforces full feasibility with fast local convergence.
Equality and inequality constraints are handled through an augmented Lagrangian formulation, and a pseudo-Huber homotopy is used to improve robustness for fuel-optimal objectives.
The solver is evaluated on benchmark transfers in Sun-centred, Earth--Moon, and Earth-centred dynamical environments.
Across these cases, the most robust configuration (iLQRDyn) converged systematically and reduced run time by 70\% (Sun-centred), 51--88\% (Earth--Moon), and 41--55\% (Earth-centred) relative to the corresponding baseline.
When convergent, the \gls*{DDP}-based variants are faster still. Overall, the results show that \gls*{DA}-based acceleration can substantially improve practical efficiency while retaining robust convergence behaviour on the tested benchmark set.},
keywords = {low-thrust trajectory optimisation \sep differential algebra \sep differential dynamic programming \sep iLQR \sep augmented Lagrangian methods}
}

\begin{document}

\maketitle

\section{Introduction}
\glsresetall

Newly developed space missions aim to explore complex environments with non-linear dynamics, such as the Earth-Moon system \citep{ZhengEtAl_2008_CLEPPaF,SmithEtAl_2020_TAPaOoNAtRHttM}. 
However, optimisation solvers may struggle in these regimes, being both time-consuming and highly sensitive to parameter settings \citep{CalebEtAl_2024_FGEaMMfBITwPO}.
Current non-linear trajectory optimisation solvers are computationally intensive\tr{. While faster methods exist, they typically rely on linear approximations of the dynamics, which can be inaccurate \citep{BooneMcMahon_2021_OGUHOSTT}.}

\tr{\Gls*{DDP}, first introduced by \citet{Mayne_1966_ASOGMfDOToNLDTS}, is a trajectory optimisation method that has been applied in mission design, including NASA’s Dawn \tg{mission} \citep{RaymanEtAl_2006_DaMiDfEoMBAVaC}\tg{(using the solver Mystic \citep{Whiffen_2006_MIotSDOCAfHFLTTD})} and \tg{the} Psyche mission \citep{OhEtAl_2017_DotPMfNDP}. Its robustness to poor initial guesses and adaptability make it attractive for a range of optimal control problems. Many extensions of the original framework have since been proposed \citep{LantoineRussell_2012_AHDDPAfCOCPP1T,OzakiEtAl_2018_SDDPwUTfLTTD}.}
\tp{Numerous solvers implement constraints using an \gls*{AUL} formulation. This approach adds the constraints to the cost function via a dual state and a penalty term that enforces the constraints while minimising the objective \citep{LinArora_1991_DDPTfCOCP1TD,Ruxton_1993_DDPAtCOCPwSVIC}. \citet{HowellEtAl_2019_AaFSfCTO} employed this formulation followed by a Newton‑based polishing phase to achieve high‑precision feasibility with shorter run times.}

\tr{\gls*{DDP} faces several challenges when applied to non-linear optimal control problems: efficiently computing second-order derivatives, handling state and control constraints, and updating the state through repeated evaluations of the system dynamics. Existing approaches often address these challenges in isolation. Constrained variants such as \gls*{HDDP}, the interior-point \gls*{DDP} method by \citet{PavlovEtAl_2021_IPDDP}, or \citet{XieEtAl_2017_DDPwNC}'s active-set strategy all assume access to second-order derivatives which is usually exact and costly or efficient and inexact approximation.
To address this, automatic differentiation techniques have been developed. For instance, \citet{NgangaWensing_2021_ASODDPfRBS} accelerate \gls*{DDP} by using reverse-mode automatic differentiation to compute second-order terms in the dynamics efficiently, avoiding the need to explicitly form large derivative tensors.
\citet{MaestriniEtAl_2018_HDDPAfLTTDUEHOTM} introduced the use of Taylor polynomials in \gls*{DDP} for derivative computation and control updates via polynomial inversion \citep{Berz_1999_MMMiPBP,ValliEtAl_2013_NMoUiCM}. However, their approach increases computational complexity by including the Lagrange multipliers as polynomial variables, leading to costly algebraic operations \citep{VasileEtAl_2019_SPiDSwGPAaICC}. In addition, even though this approach enables the use of third- or fourth-order \gls*{DDP}, their results show that the additional computations primarily lead to slower run times without significant benefits. Therefore, second-order \gls*{DDP} appears to offer a favourable compromise.
Separately, \citet{BooneMcMahon_2025_RSTOUSTTaDDP} show that a major portion of \gls*{DDP}'s runtime is spent propagating the state through non-linear dynamics.
To address this, high-order Taylor expansions have been used to approximate the dynamics and accelerate trajectory updates \citep{ArmellinEtAl_2010_ACECUDAtCoA}, but this technique typically requires an initial guess close to the optimal solution to be effective.}

\tp{In this work, we propose the \gls*{DADDy} solver, a unified and publicly available framework that addresses all three challenges.  
The framework handles fuel-optimal trajectory optimisation with equality and inequality constraints through an Augmented Lagrangian formulation, implemented in an outer loop to avoid additional computational overhead.  
The solver combines two complementary components:  
\begin{enumerate}
    \item a \gls*{DDP} solver, which generates an initial almost-feasible trajectory without requiring a good initial guess, and  
    \item a Newton solver \citep{HowellEtAl_2019_AaFSfCTO}, which polishes the solution to full feasibility while preserving optimality. 
\end{enumerate}
This structure leverages each method in its most effective domain: \gls*{DDP} provides robustness in the presence of poor initial guesses, while the Newton solver achieves machine-precision feasibility with quadratic convergence to reduce run time.  
Furthermore, we reduce the computational complexity of the Newton solver when compared to \citet{HowellEtAl_2019_AaFSfCTO} by one order by exploiting the structure of optimal control problems.  
Finally, both solvers are substantially accelerated through the use of Taylor polynomial expansions, enabling efficient derivative computation and accurate dynamics approximations. The end result is a \gls*{DDP} solver that achieves significant run time reductions on a wide variety of benchmark astrodynamics problems when compared to existing state-of-the-art \gls*{DDP} algorithm}

After introducing \tg{constrained} \gls*{DDP} and high-order Taylor polynomials in \Sec{sec:background}, \Sec{sec:methodology} details the proposed methodology.
\tp{\Sec{sec:application} first validates the solver and investigates parameter tuning, before applying the proposed methods to a range of test cases from the literature for validation and comparison with the state of the art.}
Finally, \Sec{sec:conclusions} presents the conclusions.
\section{Background} \label{sec:background}

\subsection{Constrained differential dynamic programming}

\subsubsection{Differential dynamic programming}

\Gls*{DDP} tackles optimisation problems of $N$ stages with dynamics:
\begin{equation*}
\bm{x}_{k+1}=\bm{f}\left(\bm{x}_k, \bm{u}_k\right), \quad k\in[0, N-1]
\end{equation*}
where $\bm{x}_k\in\mathbb{R}^{N_x}$ is the state vector, $\bm{u}_k\in\mathbb{R}^{N_u}$ is the control input, and $\bm{f}$ denotes the system dynamics \cite{Mayne_1966_ASOGMfDOToNLDTS}. 
The initial and target state vectors $\bm{x}_0$ and $\bm{x}_t$ of size $N_x$ are given, and the cost function to minimise is:
\begin{equation}
    J\left(\bm{U}\right) = \sum_{k=0}^{N-1} {{\ell}\left(\bm{x}_k, \bm{u}_k\right)} + {\phi}\left(\bm{x}_N, \bm{x}_t\right),
\end{equation}
where $\bm{U}=\left(\bm{u}_0, \bm{u}_1, \dotsc, \bm{u}_{N-1}\right)$ is the vector of controls of size $N N_u$, $\bm{x}_N$ is the final propagated state vector, ${\ell}$ is the stage cost, and ${\phi}$ is the terminal cost.

Let us define some notations:
\begin{enumerate}
    \item A subscript $k$ is used to indicate that dynamics, constraints, or cost functions are evaluated at the step $k$, for instance $\bm{f}_k=\bm{f}\left(\bm{x}_k, \bm{u}_k\right)$ and ${\phi}_k = {\phi}\left(\bm{x}_k,\bm{x}_t\right)$.
    \item Partial derivatives of a function $h$ with respect to a given variable $y$ are written as: $h_y=\frac{\partial h}{\partial y}$. Second{-}order derivatives with respect to $y$ then $z$: $\frac{\partial^2 h}{\partial z\partial y}$, are written $h_{yz}$.
    \item A comma between subscripts indicates the step, and then the derivative index: $\frac{\partial^2\bm{f}}{\partial \bm{x}^2}\left(\bm{x}_k, \bm{u}_k\right)=\bm{f}_{k, \bm{x}\bm{x}}$.
    \item The $i$\textendash{th} component of a vector $\bm{y}$ is written $y^i$, and the term at the $i$\textendash{th} row and $j$\textendash{th} column of a matrix $\bm{y}$ is written $y^{i,j}$.
\end{enumerate}
Let $\bm{X}=\left(\bm{x}_0, \bm{x}_1, \dotsc, \bm{x}_{N}\right)$ and define $V_k$ as the minimised cost-to-go using Bellman's principle of optimality:
\begin{equation}
    V_k\left(\bm{x}_k\right) = \min_{\bm{u}_k}\left[ {\ell}\left(\bm{x}_k, \bm{u}_k\right) + V_{k+1}\left(\bm{f}\left(\bm{x}_k, \bm{u}_k\right)\right)\right].
\end{equation}
The cost-to-go at stage $k$ is thus:
\begin{equation*}
Q_k\left(\bm{x}_k, \bm{u}_k\right)={\ell}\left(\bm{x}_k, \bm{u}_k\right) + V_{k+1}\left(\bm{f}\left(\bm{x}_k, \bm{u}_k\right)\right), \quad \text{where } V_0=\min_{\bm{U}}J.
\end{equation*}

The second-order partial derivatives of the cost-to-go are computed as:
\begin{equation}
    \label{eq:Q_derivatives}
    \begin{aligned}
        \bm{Q}_{k,\bm{x}\bm{x}} &= \bm{\ell}_{k,\bm{x}\bm{x}} + \bm{f}_{k,\bm{x}}^{\textrm{T}}\bm{V}_{k+1,\bm{x}\bm{x}}\bm{f}_{k,\bm{x}} \left(+ \sum_{i=1}^{N_x} {V_{k+1,\bm{x}}^i \left(f^i_k\right)_{\bm{x}\bm{x}}} \right),\\
        \bm{Q}_{k,\bm{x}\bm{u}} &= \bm{\ell}_{k,\bm{x}\bm{u}} + \bm{f}_{k,\bm{x}}^{\textrm{T}}\bm{V}_{k+1,\bm{x}\bm{x}}\bm{f}_{k,\bm{u}} \left(+ \sum_{i=1}^{N_x} {V_{k+1,\bm{x}}^i \left(f^i_k\right)_{\bm{x}\bm{u}}} \right),\\
        \bm{Q}_{k,\bm{u}\bm{u}} &= \bm{\ell}_{k,\bm{u}\bm{u}} + \bm{f}_{k,\bm{u}}^{\textrm{T}}\bm{V}_{k+1,\bm{x}\bm{x}}\bm{f}_{k,\bm{u}} \left(+ \sum_{i=1}^{N_x} {V_{k+1,\bm{x}}^i \left(f^i_k\right)_{\bm{u}\bm{u}}} \right).\\
    \end{aligned}
\end{equation}
The parenthetical terms represent second-order dynamics derivatives. These are neglected by \gls*{iLQR} solvers, making their cost-to-go derivatives approximate, but are included in \gls*{DDP} solvers for exactness \citep{Mayne_1966_ASOGMfDOToNLDTS, HowellEtAl_2019_AaFSfCTO}. 

This inclusion enables \gls*{DDP} to achieve quadratic convergence compared to \gls*{iLQR}'s super-linear convergence \citep{NgangaWensing_2021_ASODDPfRBS}. However, this improved convergence may reduce robustness to poor initial guesses unless proper regularisation is employed.
The first-order terms $\bm{Q}_{k,\bm{x}}$ and $\bm{Q}_{k,\bm{u}}$ are given in \citet{Mayne_1966_ASOGMfDOToNLDTS}. Computations begin with terminal conditions $\bm{V}_{N,\bm{x}} = \bm{{\phi}}_{{\bm{x}}}$ and $\bm{V}_{N,\bm{x}\bm{x}} = \bm{{\phi}}_{{\bm{x}\bm{x}}}$, and proceed sequentially backward in time, precluding parallelisation.

The \gls*{DDP} algorithm consists of two phases:
\begin{enumerate}[leftmargin=*]
\item \textbf{Backward sweep}: Compute control corrections $\bm{a}_k$ (feed-forward terms) and $\bm{K}_k$ (feedback gains) by minimising the quadratic approximation of $Q_k$ (Eq.~\ref{eq:Q_derivatives}) from $k=N-1$ to $0$ \citep{Mayne_1966_ASOGMfDOToNLDTS}.
\item \textbf{Forward pass}: Update trajectory using the computed corrections \citep{Mayne_1966_ASOGMfDOToNLDTS}.
\end{enumerate}
The backward sweep is detailed in \Alg{alg:backward_sweep}, where parenthetical instructions apply specifically to \gls*{DDP}. To ensure positive definitiveness of $\bm{Q}_{k,\bm{u}\bm{u}}$, a regularisation term is added \citep{HowellEtAl_2019_AaFSfCTO}; this is omitted from subsequent algorithms for clarity.
\begin{algorithm}
    \caption{Backward sweep}
    \label{alg:backward_sweep}
    \begin{algorithmic}
        \State \textbf{Input:} $\bm{X}$, $\bm{U}$, $\bm{x}_t$, $\bm{f}$, ${\ell}$, ${\phi}$
        \State $k \gets N-1$
        \State Compute $\bm{V}_{N,\bm{x}}$ and $\bm{V}_{N,\bm{x}\bm{x}}$
        \While{$k\geq0$}
             \State Compute $\bm{f}_{k,\bm{x}}$, $\bm{f}_{k,\bm{u}}${(, $\bm{f}_{k,\bm{x}\bm{x}}$, $\bm{f}_{k,\bm{x}\bm{u}}$, $\bm{f}_{k,\bm{u}\bm{u}}$)}
             \State Compute $\bm{{\ell}}_{k,\bm{x}}$, $\bm{{\ell}}_{k,\bm{u}}$, $\bm{{\ell}}_{k,\bm{x}\bm{x}}$, $\bm{{\ell}}_{k,\bm{x}\bm{u}}$, $\bm{{\ell}}_{k,\bm{u}\bm{u}}$
             \State Compute $\bm{Q}_{k,\bm{x}}$, $\bm{Q}_{k,\bm{u}}$, $\bm{Q}_{k,\bm{x}\bm{x}}$, $\bm{Q}_{k,\bm{x}\bm{u}}$, $\bm{Q}_{k,\bm{u}\bm{u}}$ {(using the second{-}order derivatives of $\bm{f}_{k}$)}
             \State $\bm{a}_{k},\ \bm{K}_{k}\gets-\bm{Q}_{k, \bm{u}\bm{u}}^{-1}\bm{Q}_{k,\bm{u}}^{{\textrm{T}}},\ -\bm{Q}_{k, \bm{u}\bm{u}}^{-1}\bm{Q}_{k,\bm{x}\bm{u}}^{{\textrm{T}}}$
             \State Compute $\bm{V}_{k,\bm{x}}$ and $\bm{V}_{k,\bm{x}\bm{x}}$
             \State $k \gets k - 1$
        \EndWhile
        \State \textbf{Return:} $\bm{A}$, $\bm{K}$
    \end{algorithmic}
\end{algorithm}

In the forward pass, the trajectory is updated using the computed corrections $\bm{A}$ and $\bm{K}$, producing new states $\bm{x}^*_k$ and controls $\bm{u}^*_k$ (with $\bm{x}^*_0=\bm{x}_0$). These are computed iteratively as:
\begin{equation}
    \label{eq:forward_pass}
    \begin{aligned}
        \delta\bm{x}_k^* & = \bm{x}_k^* - \bm{x}_k, \\
        \bm{u}_k^* & = \bm{u}_k + \delta\bm{u}_k^* = \bm{u}_k + \bm{a}_k + \bm{K}_k\delta\bm{x}_k^{*}, \\
        \bm{x}_{k+1}^* & =  \bm{f}\left(\bm{x}_k^*, \bm{u}_k^*\right). \\
    \end{aligned}
\end{equation}
The new cost $J^{*}$ is computed during this phase. The forward pass is detailed in \Alg{alg:forward_pass}, and the complete \gls*{DDP} algorithm in \Alg{alg:DDP}. In practice, line-search scaling $\alpha \in (0,1]$ on feedforward terms $\bm{a}_k$ improves convergence \citep{HowellEtAl_2019_AaFSfCTO}, but is omitted here for clarity.
\begin{algorithm}[h]
    \caption{Forward pass}
    \label{alg:forward_pass}
    \begin{algorithmic}
        \State \textbf{Input:} $\bm{X}$, $\bm{U}$, $\bm{x}_t$, $\bm{A}$, $\bm{K}$, $\bm{f}$, ${\ell}$, ${\phi}$
        \State $k \gets 0$
        \State $\bm{x}_0^*,\ J^* \gets \bm{x}_0,\ 0$    
        \While{$k\leq N-1$}
             \State $\delta\bm{x}_{k}^*\gets\bm{x}_{k}^*-\bm{x}_{k}$
             \State $\bm{u}_k^* \gets \bm{u}_k + \bm{a}_k + \bm{K}_k\delta\bm{x}_k^{*}$
             \State $\bm{x}_{k+1}^* \gets \bm{f}\left(\bm{x}_k^*, \bm{u}_k^*\right)$
             \State $J^*\gets J^* + {\ell}\left(\bm{x}_k^*, \bm{u}_k^*\right)$
             \State $k \gets k + 1$
        \EndWhile
        \State $J^*\gets J^* + {\phi}\left(\bm{x}_N^*, \bm{x}_t\right)$
        \State \textbf{Return:} $J^*$, $\bm{X}^*$, $\bm{U}^*$
    \end{algorithmic}
\end{algorithm}
\begin{algorithm}[h]
    \caption{\Gls*{DDP} solver}
    \label{alg:DDP}
    \begin{algorithmic}
        \State \textbf{Input:} $\varepsilon_{\textrm{DDP}}$, $\bm{x}_0$, $\bm{x}_t$, $\bm{U}_0$, $\bm{f}$, ${\ell}$, ${\phi}$
        \State $J, \ \bm{U}  \gets -\infty, \ \bm{U}_0$
        \State Compute $\bm{X}$ and $J^*$
        \While{$J^* < J\land |J-J^*|>\varepsilon_{\textrm{DDP}}$}
            \State $J\gets J^*$
            \State $\bm{A}, \ \bm{K} \gets$ BackwardSweep$\left(\bm{X}, \bm{U}, \bm{x}_t, \bm{f}, {\ell}, {\phi}\right)$
            \State $J^*,\ \bm{X}^*,\ \bm{U}^* \gets$ ForwardPass$\left(\bm{X}, \bm{U}, \bm{x}_t, \bm{A}, \bm{K}, \bm{f}, {\ell}, {\phi}\right)$    
            \State $\bm{X},\ \bm{U} \gets \bm{X}^*,\ \bm{U}^*$     
        \EndWhile
        \State \textbf{Return:}  $J^*$, $\bm{X}$, $\bm{U}$
    \end{algorithmic}
\end{algorithm}

\subsubsection{Augmented Lagrangian formulation}

{The original \gls*{DDP} algorithm presented in \citet{Mayne_1966_ASOGMfDOToNLDTS} does not include a formulation to include constraints. This work implements them using the strategy developed by \citet{HowellEtAl_2019_AaFSfCTO}. They can handle path constraints:}
\begin{equation}
    \begin{aligned}
        \bm{g}_{{\textrm{ineq}}}\left(\bm{x}_k, \bm{u}_k\right)  & \preceq \bm{0}, \\
        \bm{g}_{{\textrm{eq}}}\left(\bm{x}_k, \bm{u}_k\right)  & = \bm{0}, \\
    \end{aligned}
\end{equation}
where $\bm{g}_{{\textrm{ineq}}}$ are the inequality path constraints of size $N_{{\textrm{ineq}}}$, $\bm{g}_{{\textrm{eq}}}$ are the equality path constraints of size $N_{{\textrm{eq}}}$, and $\bm{y}\preceq\bm{0}$ means that every component of $\bm{y}$ is negative. They also implement terminal constraints:
\begin{equation}
    \begin{aligned}
        \bm{g}_{{\textrm{tineq}}}\left(\bm{x}_N, \bm{x}_t\right)  & \preceq \bm{0}, \\
        \bm{g}_{{\textrm{teq}}}\left(\bm{x}_N, \bm{x}_t\right)  & = \bm{0}, \\
    \end{aligned}
\end{equation}
where $\bm{g}_{{\textrm{tineq}}}$ are the terminal inequality constraints of size $N_{{\textrm{tineq}}}$, and $\bm{g}_{{\textrm{teq}}}$ are the terminal equality constraints of size $N_{{\textrm{teq}}}$. Let us define: $\bm{g} =\left[\bm{g}_{{\textrm{ineq}}}^{{\textrm{T}}}, \bm{g}_{{\textrm{eq}}}^{{\textrm{T}}}\right]^{{\textrm{T}}}$, $\bm{g}_t = \left[\bm{g}_{{\textrm{tineq}}}^{{\textrm{T}}}, \bm{g}_{{\textrm{teq}}}^{{\textrm{T}}}\right]^{{\textrm{T}}}$, and $\bm{g}_N =\left[\bm{g}_{{\textrm{tineq}}}\left(\bm{x}_N, \bm{x}_t\right)^{{\textrm{T}}}, \bm{g}_{{\textrm{teq}}}\left(\bm{x}_N, \bm{x}_t\right)^{{\textrm{T}}}\right]=\bm{g}_t\left(\bm{x}_N, \bm{x}_t\right)$.
These constraints are handled using the \gls*{AUL} formulation.
{This approach relies on the fact that a constrained optimisation problem can be equivalently reformulated as an unconstrained one, where constraint satisfaction is {handled} by adding terms to the cost functions. These additional terms combine dual variables (Lagrange multipliers) with quadratic penalties on the constraint violations.
As a result, the optimiser is guided toward feasibility without requiring explicit constraint enforcement at each iteration.} 
{Following the methodology developed in \citet{HowellEtAl_2019_AaFSfCTO}, the functions $\ell$ and $\phi$ are augmented as $\tilde{\ell}$ and $\tilde{\phi}$ using dual states and penalty vectors denoted respectively $\bm{\Lambda} = \left(\bm{\lambda}_0, \bm{\lambda}_1, \dotsc, \bm{\lambda}_N\right)$ and $\bm{M} = \left(\bm{\mu}_0, \bm{\mu}_1, \dotsc, \bm{\mu}_N\right)$.}
{Let us define the vector of constraints} $\bm{G}=\left(\bm{g}_0, \bm{g}_1, \dotsc, \bm{g}_N\right)$.
The augmented problem is solved using \gls*{DDP} until the maximum constraints violation $g_{\max}{=\max \bm{G}}$ is below the target constraint satisfaction $\varepsilon_{\textrm{AUL}}>0$.

\subsection{Solution polishing using a Newton method}
\label{sec:background_newton}
Solving fuel-optimal problems with \gls*{DDP} typically requires final iterations to fine-tune controls for constraint satisfaction. This refinement often consumes over half the total \gls*{DDP} iterations, making it computationally expensive.

\citet{HowellEtAl_2019_AaFSfCTO} propose a two-stage strategy: (1) terminate the augmented Lagrangian \gls*{DDP} solver at low-accuracy constraint satisfaction level $\varepsilon_{\textrm{AUL}}$ once cost converges; (2) refine via a Newton-based polishing step to achieve high-precision feasibility $\varepsilon_{\textrm{N}} \ll \varepsilon_{\textrm{AUL}}$ while maintaining cost accuracy. Since Newton exhibits quadratic convergence near optimality, this avoids the expensive iterative refinement of standard \gls*{DDP}.

The method concatenates all active constraints (including continuity) into a vector $\bm{\Gamma}$ of size $N_\Gamma \leq N\cdot (N_{{\textrm{ineq}}} + N_{{\textrm{eq}}} + N_x) + N_{{\textrm{tineq}}} + N_{{\textrm{teq}}}$: 
\begin{equation}
    \bm{h}_k=\bm{h}\left(\bm{x}_k, \bm{u}_k, \bm{x}_{k+1}\right)=\bm{x}_{k+1} -  \bm{f}\left(\bm{x}_k, \bm{u}_k\right)=\bm{0}.
\end{equation}
Free variables (all states and controls except $\bm{x}_0$) are concatenated into $\bm{Y}$ of size $N_Y=N\cdot(N_x+N_u)$. The constraints are linearised as:
\begin{equation*}
\bm{\Gamma}(\bm{Y} + \delta\bm{Y}) \approx \bm{\Delta}\delta\bm{Y} + \bm{d},
\end{equation*}
where $\bm{\Delta}$ and $\bm{d}$ are the constraint Jacobian and value at $\bm{Y}$. The Newton correction $\delta\bm{Y}^* = -\bm{\Delta}^{+} \bm{d}$ (using pseudo-inverse $\bm{\Delta}^{+}$) solves the non-square system $\bm{\Delta}  \bm{z} + \bm{d} = \bm{0}$ to yield refined estimate $\bm{Y}^*$.
The entire process is presented in \Alg{alg:pn_method}.
\begin{algorithm}
    \caption{Newton solver}
    \label{alg:pn_method}
    \begin{algorithmic}
        \State \textbf{Input:} $\varepsilon_{\textrm{N}}$, $\bm{X}_0$, $\bm{U}_0$, $\bm{x}_t$, $\bm{f}$, $\bm{g}$, $\bm{g}_t$
        \State $\gamma, \ d_{\max} \gets 0.5, \ +\infty$  
        \State Compute $\bm{Y}$
        \State Retrieve $\bm{d}$
        \While{$d_{\max}>\varepsilon_{\textrm{N}}$}
            \State Compute $\bm{\Delta}$ at $\bm{Y}$
            \State $r\gets +\infty$ 
            \While{$d_{\max}>\varepsilon_{\textrm{N}}\land r>\varepsilon_{\textrm{CV}}$} \Comment{Reuse $\bm{\Delta}$}
                \State $d_{\max}^*,\ \alpha\gets +\infty,\ 1$ 
                \While{$d_{\max}^*>d_{\max}$} \Comment{Linesearch}
                    \State $\bm{Y}^*\gets \bm{Y}-\alpha{\bm{\Delta}^{+}}\bm{d}$
                    \State $\bm{d}^*\gets\bm{\Gamma}\left(\bm{Y}^*\right)$
                    \State $d_{\max}^*,\ \alpha\gets \max |\bm{d}^*|, \ \gamma\alpha$
                \EndWhile
                \State $r, \ d_{\max}\gets \log d_{\max}^*/\log d_{\max}, \ d_{\max}^*$
                \State $\bm{d}, \ \bm{Y}\gets\bm{d}^*,\ \bm{Y}^*$
            \EndWhile
        \EndWhile
        \State \textbf{Return:}  $\bm{X}$, $\bm{U}$, $d_{\max}$
    \end{algorithmic}
\end{algorithm}
The parameter $r$ represents the convergence rate per iteration, typically with threshold $\varepsilon_{\textrm{CV}}=\num{1.1}$ \citep{HowellEtAl_2019_AaFSfCTO}. This allows reusing the expensive matrix $\bm{\Delta}$ until convergence degrades, then recomputes it for acceleration. The result is a solution satisfying both optimality ($\varepsilon_{\textrm{DDP}}$) and feasibility ($\varepsilon_{\textrm{N}}$) tolerances.

\subsection{Differential Algebra}
\Gls*{DA} is the dedicated framework to manipulate high-order Taylor expansions \cite{Berz_1999_MMMiPBP}. {Though it was initially designed to efficiently perform high-order automatic differentiation, many other applications followed such as uncertainty propagation or surrogate modelling, as presented in \citet{Wittig_2016_AItDAatDAMR}.
This framework enables the representation of a given sufficiently differentiable function $\bm{h}$ of $v$ variables with its Taylor expansion $\mathcal{P}_h$. 
Furthermore, algebraic and functional operations, including polynomial inversion, are well-defined on the set of all Taylor polynomials \cite{Berz_1992_HOCaNFAoRS}. }

{In order to perform automatic differentiation using \gls*{DA}, t}he variable is defined as a vector of polynomials: $\mathcal{P}_{v} = \bm{v} + \delta\bm{v}${, w}here $\bm{v}$ is the constant part and $\delta\bm{v}$ is a perturbation vector.
Evaluating a function $\bm{h}$ at point $\mathcal{P}_{v}$ returns a Taylor expansion $\mathcal{P}_{h}$ of $\bm{h}$ around $\bm{v}$. The derivatives of $\bm{h}$ at $\bm{v}$ can be retrieved {from} the coefficients of $\mathcal{P}_{h}$ {at negligible computational cost.}

Throughout the paper, the constant part of a polynomial $\mathcal{P}_h$ is denoted $\overline{\mathcal{P}_{h}}${, and} a small, unset, perturbation of any given quantity $y$ is denoted $\delta y$.
{For all $\varepsilon>0$, {we can define a convergence radius $R_{\varepsilon}$} such that, given a Taylor expansion $\mathcal{P}_{h}$ of $\bm{h}$ around point $\bm{v}$ we have:}
\begin{equation}
    \label{eq:convergence_radius}
    \|\delta\bm{v}\|_2 \leq R_{\varepsilon} \implies \left\|\mathcal{P}_{h}\left(\delta\bm{v}\right) - \bm{h}\left(\bm{v}+\delta\bm{v}\right)\right\|_2 \leq \varepsilon.
\end{equation}
{I}n other words, $\mathcal{P}_{h}$ is a local approximation of $\bm{h}$ around the point $\bm{v}$. 
{The method used in this work to compute the convergence radius follows the approach proposed by \citet{WittigEtAl_2015_PoLUSiODbADS}, which is based on the exponential decay of the Taylor coefficients with increasing order.}
This feature proves useful when evaluations of the function $h$ are costly{, for instance in state propagation \citep{ArmellinEtAl_2010_ACECUDAtCoA}, complex orbit propagations \citep{CalebEtAl_2022_SSMwTDAwAtBCOaM,CalebEtAl_2023_DAMAtCAGaBDAtPFotES}, or uncertainty propagation \citep{WittigEtAl_2015_PoLUSiODbADS, LosaccoEtAl_2024_LOADSAfNUM}.} {In such cases, the use of a polynomial approximation can significantly reduce run time if a significant number of function evaluations are needed, as in trajectory optimisation \citep{CalebEtAl_2024_FGEaMMfBITwPO} or Monte-Carlo analyses \citep{ValliEtAl_2013_NMoUiCM}}. {In these studies, \gls*{DA} is used for both automatic differentiation and to obtain accurate approximations of the periodic orbits of the \gls*{CR3BP}.}

\section{Methodology} \label{sec:methodology}

This section \tp{first details the use of} high-order Taylor polynomials in \gls*{DDP} in the \gls*{DADDy} solver \tp{for both automatic differentiation and to approximate the repeated evaluations of the dynamics}.
Then, a method to perform fuel-optimal optimisation is shown, followed by an enhanced solution polishing method.

\subsection{\Gls*{DA}-based \gls*{DDP}}
\label{sec:DA_DDP}
\subsubsection{Automatic differentiation}
\tp{The first use of \gls*{DA} is the \gls*{DADDy} algorithm is for automatic differentiation.} States and controls are defined as vectors of polynomials: $\mathcal{P}_{x} = \bm{x} + \delta\bm{x}$ and $\mathcal{P}_{u} = \bm{u} + \delta\bm{u}$, where $\bm{x}$ (respectively $\bm{u}$) is a computed state vector (control) and $\delta\bm{x}$ ($\delta\bm{u}$) is a perturbation vector of size $N_x$ ($N_u$).
Therefore, evaluating any function $\bm{h}$ at point $\left(\mathcal{P}_{x}, \mathcal{P}_{u}\right)$ at order $2$ returns a second-order Taylor expansion $\mathcal{P}_{h}$ of $\bm{h}$ at $\left(\bm{x}, \bm{u}\right)$.
Then, its first-order derivatives $\bm{h}_{\bm{x}}$ and $\bm{h}_{\bm{u}}$, and second-order derivatives $\bm{h}_{\bm{x}\bm{x}}$, $\bm{h}_{\bm{x}\bm{u}}$, and $\bm{h}_{\bm{u}\bm{u}}$ are retrieved \tr{from} the coefficients of $\mathcal{P}_{h}$. 
Thus, the derivatives of $\bm{f}$, $\tilde{{\ell}}$, and $\tilde{\tr{\phi}}$ can be retrieved without implementing the derivatives by hand or being limited to linear and quadratic problems \cite{HowellEtAl_2019_AaFSfCTO}.
The previous considerations result in a polynomial-based forward pass, shown in \Alg{alg:forward_pass_AD}. It allows for automatic differentiation of the dynamics, stage cost, and terminal cost functions, as they are evaluated using the \gls*{DA} variables $\delta \bm{x}$, and $\delta \bm{u}$. As a shorthand, the following notation is used: $\mathcal{P}_{{\ell}}=\left(\mathcal{P}_{{\ell},0}, \mathcal{P}_{{\ell},1},\dotsc, \mathcal{P}_{{\ell}, N-1}\right)$, similarly for $\mathcal{P}_{f}$.
\begin{algorithm}
    \caption{Forward pass with automatic differentiation}\label{alg:forward_pass_AD}
    \begin{algorithmic}
        \State \textbf{Input:} $\bm{X}$, $\bm{U}$, $\bm{x}_t$, $\bm{A}$, $\bm{K}$, $\bm{f}$, ${\ell}$, $\tr{\phi}$
        \State $k \gets 0$
        \State $\bm{x}_k^*,\ J^* \gets \bm{x}_0,\ 0$
        \While{$k\leq N-1$}
            \State $\delta\bm{x}_{k}^*\gets\bm{x}_{k}^*-\bm{x}_{k}$
            \State $\bm{u}_k^* \gets \bm{u}_k + \bm{a}_k + \bm{K}_k\delta\bm{x}_k^{*}$
            \State $\mathcal{P}_{f,k},\ \mathcal{P}_{{\ell},k} \gets \bm{f}\left(\bm{x}_k^* + \delta\bm{x}, \bm{u}_k^* + \delta\bm{u}\right), \  {\ell}\left(\bm{x}_k^* + \delta\bm{x}, \bm{u}_k^* + \delta\bm{u}\right)$
            \State $\bm{x}_{k+1}^*,\ J^* \gets \overline{\mathcal{P}_{f,k}},\ J^* + \overline{\mathcal{P}_{{\ell},k}}$            
            \State $k \gets k + 1$
        \EndWhile
        \State $\mathcal{P}_{\tr{\phi}}\gets \tr{\phi}\left(\bm{x}_N^* + \delta\bm{x}, \bm{x}_t\right)$
        \State $J^*\gets J^* + \overline{\mathcal{P}_{\tr{\phi}}}$
        \State \textbf{Return:} $J^*$, $\bm{X}^*$, $\bm{U}^*$, $\mathcal{P}_{f}$, $\mathcal{P}_{{\ell}}$, $\mathcal{P}_{\tr{\phi}}$
    \end{algorithmic}
\end{algorithm}
Consequently, a backward sweep using the automatic differentiation of \Alg{alg:forward_pass_AD} is shown in \Alg{alg:backward_sweep_AD}\tr{, where the instructions between parentheses refer specifically to \gls*{DDP}.}
\tr{Indeed, as mentioned in \Sec{sec:background}, \gls*{DDP} requires second-order derivatives of the dynamics while \gls*{iLQR} solvers do not.
These derivatives are already computed as part of the dynamics approximation and can be added without requiring additional computations or numerical integrations.}
\begin{algorithm}
    \caption{Backward sweep with automatic differentiation}\label{alg:backward_sweep_AD}
    \begin{algorithmic}
        \State \textbf{Input:} $\bm{X}$, $\bm{U}$, $\bm{x}_t$, $\mathcal{P}_{f}$, $\mathcal{P}_{{\ell}}$, $\mathcal{P}_{\tr{\phi}}$
        \State Retrieve $\bm{V}_{N,\bm{x}}$ and $\bm{V}_{N,\bm{x}\bm{x}}$ from $\mathcal{P}_{\tr{\phi}}$
        \State $k \gets N-1$
        \While{$k\geq0$}
             \State Retrieve $\bm{f}_{k,\bm{x}}$, $\bm{f}_{k,\bm{u}}$\tr{(, $\bm{f}_{k,\bm{x}\bm{x}}$, $\bm{f}_{k,\bm{x}\bm{u}}$, $\bm{f}_{k,\bm{u}\bm{u}}$)} from $\mathcal{P}_{f,k}$
             \State Retrieve $\bm{{\ell}}_{k,\bm{x}}$, $\bm{{\ell}}_{k,\bm{u}}$, $\bm{{\ell}}_{k,\bm{x}\bm{x}}$, $\bm{{\ell}}_{k,\bm{x}\bm{u}}$, $\bm{{\ell}}_{k,\bm{u}\bm{u}}$ from $\mathcal{P}_{{\ell},k}$
             \State Compute $\bm{Q}_{k,\bm{x}}$, $\bm{Q}_{k,\bm{u}}$, $\bm{Q}_{k,\bm{x}\bm{x}}$, $\bm{Q}_{k,\bm{x}\bm{u}}$, $\bm{Q}_{k,\bm{u}\bm{u}}$ \tr{(using the second\tr{-}order derivatives of $\bm{f}_{k}$)}
             \State $\bm{a}_{k},\ \bm{K}_{k}\gets-\bm{Q}_{k, \bm{u}\bm{u}}^{-1}\bm{Q}_{k,\bm{u}}^{\tr{\textrm{T}}},\ -\bm{Q}_{k, \bm{u}\bm{u}}^{-1}\bm{Q}_{k,\bm{x}\bm{u}}^{\tr{\textrm{T}}}$
             \State Compute $\bm{V}_{k,\bm{x}}$ and $\bm{V}_{k,\bm{x}\bm{x}}$
             \State $k \gets k - 1$
        \EndWhile
        \State \textbf{Return:} $\bm{A}$, $\bm{K}$
    \end{algorithmic}
\end{algorithm}

\subsubsection{Enhanced forward pass: Approximation of the dynamics}
\label{sec:DA_FP}
\tp{The second use of \gls*{DA} is the \gls*{DADDy} algorithm is for efficient approximation of the dynamics.} \citet{BooneMcMahon_2025_RSTOUSTTaDDP} highlight that the forward pass is one order of magnitude slower than the backward sweep due \tp{to the costly repeated evaluations of the dynamics}. It is the reason why they use \gls*{DA} (or \gls*{STTs}) to accelerate optimisation, given an initial trajectory.
In this work, this principle is applied at every iteration of \gls*{DDP} to accelerate the evaluation of the dynamics.
The polynomial expansion of the dynamics $\mathcal{P}_{f,k}\left(\delta\bm{x},\delta\bm{u}\right)$ at stage $k$ is already computed for the automatic differentiation and can be leveraged to avoid recomputing the dynamics when the corrections $\left(\delta\bm{x}_k^*, \delta\bm{u}_k^*\right)$ computed in \Eq{eq:forward_pass} are small. They need to be within the convergence radius $R_{\varepsilon_{\textrm{DA}},k}$ of $\mathcal{P}_{f,k}$, defined in \Eq{eq:convergence_radius}, where the accuracy is $\varepsilon_{\textrm{DA}}>0$ and $\|\delta\bm{v}\|_2=\sqrt{\|\delta\bm{x}_k^*\|^2_2 + \|\delta\bm{u}_k^*\|^2_2}$.
Finally, the dynamics are now evaluated as:
\begin{equation}
    \begin{aligned}
        \bm{x}_{k+1}^* & \approx \mathcal{P}_{f,k}\left(\delta\bm{x}_k^*,\delta\bm{u}_k^*\right), \ \text{if }  \sqrt{\|\delta\bm{x}_k^*\|^2_2 + \|\delta\bm{u}_k^*\|^2_2} \leq R_{\varepsilon_{\textrm{DA}},k},\\
        \bm{x}_{k+1}^* & = \overline{\bm{f}\left(\mathcal{P}_{\bm{x}_k^*},\mathcal{P}_{\bm{u}_k^*}\right)}, \ \text{otherwise.} \\
    \end{aligned}
\end{equation}
Note that the second case leads to recomputing the dynamics and their Taylor expansions, while the first case only consists \tr{of} performing $N_x$ polynomial compositions of $N_x+N_u$ variables to obtain the dynamics and their associated derivatives.
This approach leads to the novel polynomial-based forward pass of \Alg{alg:forward_pass_DA}. It also allows for automatic differentiation and implements polynomial approximation of the dynamics using the already computed polynomial representations of $\bm{f}$.
\begin{algorithm}
    \caption{Forward pass with dynamics approximation and automatic differentiation}\label{alg:forward_pass_DA}
    \begin{algorithmic}
        \State \textbf{Input:} $\varepsilon_{\textrm{DA}}$, $\bm{X}$, $\bm{U}$, $\bm{x}_t$, $\bm{A}$, $\bm{K}$, $\mathcal{P}_{f}$, ${\ell}$, $\tr{\phi}$
        \State $k \gets 0$
        \State $J^*,\  \bm{x}_k^* \gets 0, \ \bm{x}_0$
        \While{$k\leq N-1$}
            \State $\delta\bm{x}_{k}^*\gets\bm{x}_{k}^*-\bm{x}_{k}$
            \State $\delta\bm{u}_k^* \gets \bm{a}_k + \bm{K}_k\delta\bm{x}_k^{*}$
            \State $\bm{u}_k^* \gets\bm{u}_k + \delta\bm{u}_k^* $
            \State Compute $R_{\varepsilon_{\textrm{DA}},k}$
            \If{$\sqrt{\|\delta\bm{x}_k^*\|^2_2 + \|\delta\bm{u}_k^*\|^2_2}<R_{\varepsilon_{\textrm{DA}},k}$}
                \State $\mathcal{P}_{f,k} \gets \mathcal{P}_{f,k}\left(\delta\bm{x}_k^* + \delta\bm{x}, \delta\bm{u}_k^* + \delta\bm{u}\right)$
            \Else
                \State $\mathcal{P}_{f,k} \gets \bm{f}\left(\bm{x}_k^* + \delta\bm{x}, \bm{u}_k^*+ \delta\bm{u}\right)$
            \EndIf
            \State $\mathcal{P}_{{\ell},k}\gets {\ell}\left(\bm{x}_k^* + \delta\bm{x}, \bm{u}_k^* + \delta\bm{u}\right)$
            \State $\bm{x}_{k+1}^*,\ J^* \gets \overline{\mathcal{P}_{f,k}},\ J^* + \overline{\mathcal{P}_{{\ell},k}}$            
            \State $k \gets k + 1$
        \EndWhile
        \State $\mathcal{P}_{\tr{\phi}}\gets \tr{\phi}\left(\bm{x}_N^* + \delta\bm{x}, \bm{x}_t\right)$
        \State $J^*\gets J^* + \overline{\mathcal{P}_{\tr{\phi}}}$
        \State \textbf{Return:} $J^*$, $\bm{X}^*$, $\bm{U}^*$, $\mathcal{P}_{\bm{f}}$, $\mathcal{P}_{{\ell}}$, $\mathcal{P}_{\tr{\phi}}$
    \end{algorithmic}
\end{algorithm}

\subsubsection{Enhanced Backward sweep: direct cost-to-go derivatives computations}

\tr{As shown in \citet{MaestriniEtAl_2018_HDDPAfLTTDUEHOTM}, the cost-to-go} $Q_k$ can be directly evaluated as $\mathcal{P}_{Q,k} =  {\ell}\left(\mathcal{P}_{x,k}, \mathcal{P}_{u,k}\right) + V_{k+1}\left(\bm{f}\left(\mathcal{P}_{x,k}, \mathcal{P}_{u,k}\right)\right)$.
\tp{The partial derivatives of $Q_k$ are then extracted from the coefficients of its expansion to allow us to avoid extracting the gradients and Hessians of the dynamics, the stage cost, and the terminal cost, and to subsequently retrieve the partial derivatives of $Q_k$.}
\tr{Note that the second-order derivatives of the dynamics are automatically included.
This strategy is implemented in \Alg{alg:backward_sweep_Q}.}
\begin{algorithm}
    \caption{Backward sweep with automatic differentiation and direct $Q_k$ evaluation}\label{alg:backward_sweep_Q}
    \begin{algorithmic}
        \State \textbf{Input:} $\bm{X}$, $\bm{U}$, $\bm{x}_t$, $\mathcal{P}_{f}$, $\mathcal{P}_{{\ell}}$, $\mathcal{P}_{\tr{\phi}}$
        \State $k \gets N-1$
        \State $\mathcal{P}_{V,k}\gets\mathcal{P}_{\tr{\phi}}$
        \While{$k\geq0$}
             \State $\mathcal{P}_{Q,k}\gets \mathcal{P}_{{\ell},k} + \mathcal{P}_{V,k+1}\left(\mathcal{P}_{f,k} -\bm{x}_{k+1}, \bm{0} \right) $
             \State Retrieve $\bm{Q}_{k,\bm{x}}$, $\bm{Q}_{k,\bm{u}}$, $\bm{Q}_{k,\bm{x}\bm{x}}$, $\bm{Q}_{k,\bm{x}\bm{u}}$, $\bm{Q}_{k,\bm{u}\bm{u}}$ from $\mathcal{P}_{Q,k}$
             
             \State $\bm{a}_{k},\ \bm{K}_{k}\gets-\bm{Q}_{k, \bm{u}\bm{u}}^{-1}\bm{Q}_{k,\bm{u}}^{\tr{\textrm{T}}},\ -\bm{Q}_{k, \bm{u}\bm{u}}^{-1}\bm{Q}_{k,\bm{x}\bm{u}}^{\tr{\textrm{T}}}$
             \State $\mathcal{P}_{V,k}\gets\mathcal{P}_{Q,k}\left(\delta\bm{x}, \bm{a}_{k} + \bm{K}_{k}\delta\bm{x}\right)$
             \State $k \gets k - 1$
        \EndWhile
        \State \textbf{Return:} $\bm{A}$, $\bm{K}$
    \end{algorithmic}
\end{algorithm}

\subsection{Fuel-optimal optimisation}
\label{sec:fuel_optimal}
The energy-optimal problem is smooth and consists of a stage cost function of type: ${\ell}_{\textrm{E}}(\bm{x},\bm{u}) = \frac{\bm{u}^{\tr{\textrm{T}}}\bm{u}}{2}$. However, if $\bm{u}$ is the thrust vector of a spacecraft, the stage cost to minimise to achieve minimum fuel is of type: ${\ell}_{\textrm{F}}(\bm{x},\bm{u}) =\sqrt{\bm{u}^{\tr{\textrm{T}}}\bm{u}}$, which implies the evaluation of the square root of the control. Yet, when $\bm{u}\approx\bm{0}$  the derivatives of ${\ell}_{\textrm{F}}$ diverge since the square root is non-differentiable at $\bm{u}=\bm{0}$. \tp{Gradient-based optimisation solvers, such as \gls*{DDP}, will struggle in the vicinity of the this singularity, and the \gls*{DA} framework is no longer usable.}
\tp{Therefore, in the \gls*{DADDy} algorithm, an alternative stage cost function $\ell_\textrm{H}$ is implemented based on the pseudo-Huber loss function \cite{BaraniEtAl_2024_ADLBoRDSOANwNOD}}:
\begin{equation}
    {\ell}_{\textrm{H}}(\bm{x},\bm{u}) =\sigma\left[\sqrt{\frac{\bm{u}^{\tr{\textrm{T}}}\bm{u}}{\sigma^2} + 1} - 1 \right],
\end{equation}
where $\sigma$ is a tuning parameter. Note that:
\begin{equation}
    {\ell}_{\textrm{H}}(\bm{x},\bm{u}) \sim
    \begin{aligned}
        \frac{\bm{u}^{\tr{\textrm{T}}}\bm{u}}{2\sigma}=\dfrac{{\ell}_{\textrm{E}}(\bm{x},\bm{u})}{\sigma},\ \text{if } \frac{\bm{u}^{\tr{\textrm{T}}}\bm{u}}{\sigma^2} \ll 1, \\
          \sqrt{\bm{u}^{\tr{\textrm{T}}}\bm{u}}={\ell}_{\textrm{F}}(\bm{x},\bm{u}),\ \text{if } \frac{\bm{u}^{\tr{\textrm{T}}}\bm{u}}{\sigma^2} \gg 1. \\
    \end{aligned}
\end{equation}
Therefore, the singularity no longer exists.
The strategy to converge to the fuel-optimal problem is to use an homotopy between ${\ell}_{\textrm{E}}$ and ${\ell}_{\textrm{H}}$ \citep{JiangEtAl_2012_PTfLTTOwHA}:
\begin{equation}
    \label{eq:fuel_optimal_ctg}
    {\ell}(\bm{x},\bm{u}) = \eta {\ell}_{\textrm{E}}(\bm{x},\bm{u}) + (1-\eta) {\ell}_{\textrm{H}}(\bm{x},\bm{u}).
\end{equation}
A first round of energy-optimal optimisation is performed ($\eta=1$). Then, the goal is to converge to $\eta \approx 0$ and $\sigma\approx 0$ so that: ${\ell}\sim {\ell}_{\textrm{F}}$.

\subsection{Newton method acceleration}
\tp{The Newton method described in \Alg{alg:pn_method} and \Sec{sec:background_newton} can also be accelerated using the \gls*{DA} framework. We also propose a further acceleration method relying on the solving of symmetric block tri-diagonal systems.}

\tp{\Sec{sec:background_newton} presented the structure of the Newton solver of \citet{HowellEtAl_2019_AaFSfCTO}. This method requires repeated evaluations of the dynamics to compute the continuity constraints $\bm{h}_k$. Moreover, because the constraints are subsequently linearised, obtaining the derivatives of the vector of constraints $\bm{G}$ efficiently is essential.}
\tp{To perform repeated evaluations of the dynamics and for automatic differentiation}{, we propose to construct a polynomial approximation of the constraint vector $\bm{\Gamma}$ in terms of $\delta\bm{Y}$, the variation of the vector $\bm{Y}$. The gradient of the constraints $\bm{\Delta}$ can then be retrieved directly from the coefficients of the expansion $\mathcal{P}_{\bm{\Gamma}}$. Using the same mapping, the constraints are updated after performing a correction $\delta\bm{Y}^*$ on the variables without recomputing the true dynamics. The updated constraints are $\bm{d}^*=\mathcal{P}_{\Gamma}\left(\delta\bm{Y}^*\right)$.}

{Then}, the solving of the system $\bm{\Delta} \bm{z} \tr{+ \bm{d}=\bm{0}}$ of unknown $\bm{z}$ to obtain $\delta \bm{Y}^*$ comes down to solving the system: $\bm{\Delta}\bm{\Delta}^{\tr{\textrm{T}}} \bm{z}'=-\bm{d}$ of unknown $\bm{z}'$ {with $\bm{z}=\bm{\Delta}^{\textrm{T}}\bm{z}'$}.
If $\Tilde{\bm{g}}$ (respectively $\Tilde{\bm{h}}$) is the vector of the active constraints of $\bm{g}$ ($\bm{h}$), and $\Tilde{\mathbb{I}}_{k,N_x}$ is $\mathbb{I}_{N_x}$, the identity of size $N_x\times N_x$, after removal of the rows of the inactive components of $\bm{h}_k$. Then, the matrix $\bm{\Delta}$ is:
\begin{equation}
    \bm{\Delta}=
    \begin{bmatrix}
        \Tilde{\bm{g}}_{0,\bm{u}} & & & & & & & \\
        \Tilde{\bm{h}}_{0,\bm{u}} & -\Tilde{\mathbb{I}}_{0,N_x} & & & & & & \\
        & \Tilde{\bm{g}}_{1,\bm{x}} & \Tilde{\bm{g}}_{1,\bm{u}} & & & & & \\
        & \Tilde{\bm{h}}_{1,\bm{x}} & \Tilde{\bm{h}}_{1,\bm{u}} & -\Tilde{\mathbb{I}}_{1,N_x} & & & & \\
        & & \ddots & \ddots & \ddots & & & \\
        & & & \ddots  & \ddots & \ddots & & \\
        & & & & & \Tilde{\bm{g}}_{N-1,\bm{x}} & \Tilde{\bm{g}}_{N-1,\bm{u}} & \\
        & & & & & \Tilde{\bm{h}}_{N-1,\bm{x}} & \Tilde{\bm{h}}_{N-1,\bm{u}}  & -\Tilde{\mathbb{I}}_{N-1,N_x}  \\
        & & & & & & & \Tilde{\bm{g}}_{N,\bm{x}} \\
    \end{bmatrix}.
\end{equation}
The rest is filled with $0$, and $\bm{\Delta}$ is of size $N_c N\cdot(N_x+N_u)$, with $N_c$ the total number of active constraints.
Let us define the following notations: $\bm{I}_k=\Tilde{\mathbb{I}}_{k,N_x}\Tilde{\mathbb{I}}_{k,N_x}^{\tr{\textrm{T}}}$, $\Tilde{\bm{g}}_{k,\bm{v}}^2=\Tilde{\bm{g}}_{k,\bm{v}} \Tilde{\bm{g}}_{k,\bm{v}}^{\tr{\textrm{T}}}$, and $\Tilde{\bm{h}}_{k,\bm{v}}^2=\Tilde{\bm{h}}_{k,\bm{v}} \Tilde{\bm{h}}_{k,\bm{v}}^{\tr{\textrm{T}}}$.
Note that $\bm{\Sigma}=\bm{\Delta}\bm{\Delta}^{\tr{\textrm{T}}}$ is a symmetric block tri-diagonal matrix:
\begin{equation} 
    \bm{\Sigma} =
    \footnotesize
    \begin{bmatrix} 
        \Tilde{\bm{g}}_{0,\bm{u}}^2 & \Tilde{\bm{g}}_{0,\bm{u}}\Tilde{\bm{h}}_{0,\bm{u}}^{\tr{\textrm{T}}} & & & & & \\
        \Tilde{\bm{h}}_{0,\bm{u}}\Tilde{\bm{g}}_{0,\bm{u}}^{\tr{\textrm{T}}}  & \Tilde{\bm{h}}_{0,\bm{u}}^2 + \bm{I}_0  & -\Tilde{\bm{g}}_{1,\bm{x}}^{\tr{\textrm{T}}} & -\Tilde{\bm{h}}_{1,\bm{x}}^{\tr{\textrm{T}}}  & & & \\
         & -\Tilde{\bm{g}}_{1,\bm{x}} & \Tilde{\bm{g}}_{1,\bm{x}}^2 + \Tilde{\bm{g}}_{1,\bm{u}}^2  &  \Tilde{\bm{g}}_{1,\bm{x}}\Tilde{\bm{h}}_{1,\bm{x}}^{\tr{\textrm{T}}} + \Tilde{\bm{g}}_{1,\bm{u}}\Tilde{\bm{h}}_{1,\bm{u}}^{\tr{\textrm{T}}} & & & \\
         & -\Tilde{\bm{h}}_{1,\bm{x}} & \Tilde{\bm{h}}_{1,\bm{x}}\Tilde{\bm{g}}_{1,\bm{x}}^{\tr{\textrm{T}}} + \Tilde{\bm{h}}_{1,\bm{u}}\Tilde{\bm{g}}_{1,\bm{u}}^{\tr{\textrm{T}}} & \Tilde{\bm{h}}_{1,\bm{x}}^2 + \Tilde{\bm{h}}_{1,\bm{u}}^2 +\bm{I}_1 & \ddots & \\
         & & & \ddots & \ddots & \ddots & \\
         & & & & \ddots & \Tilde{\bm{h}}_{N - 1,\bm{x}}^2 + \Tilde{\bm{h}}_{N - 1,\bm{u}}^2 +\bm{I}_{N - 1} & -\Tilde{\bm{g}}_{N,\bm{x}}^{\tr{\textrm{T}}} \\
         & & & & & -\Tilde{\bm{g}}_{N,\bm{x}} & \Tilde{\bm{g}}_{N,\bm{x}}^2  \\
    \end{bmatrix}.
\end{equation}
It can be computed analytically from the gradients of $\bm{g}$ and $\bm{h}$ to avoid computing the product of two sparse matrices. Moreover, $\bm{\Sigma}$ is positive definite, thus, has a Cholesky factorisation $\bm{\Pi}$ such that $\bm{\Sigma}=\bm{\Pi}\bm{\Pi}^{\tr{\textrm{T}}}$ and $\bm{\Pi}$ is lower-triangular \cite{Cholesky_1910_SLRNDSDL}. It allows for faster system solving and, since $\bm{\Sigma}$ is a block tri-diagonal, its Cholesky factorisation can be retrieved faster than for most matrices of similar size \cite{CaoEtAl_2002_PCFoaBTM}.
\Alg{alg:pn_method_DA} presents these modifications.
\begin{algorithm}
    \caption{Modified Newton solver}
    \label{alg:pn_method_DA}
    \begin{algorithmic}
        \State \textbf{Input:} $\varepsilon_{\textrm{N}}$, $\bm{X}_0$, $\bm{U}_0$, $\bm{x}_t$, $\bm{f}$, $\bm{g}$, $\bm{g}_t$
        \State Compute $\bm{Y}$
        \State $d_{\max},\ \gamma \gets +\infty,\ 0.5$   
        \While{$d_{\max}>\varepsilon_{\textrm{N}}$}
            \State $\mathcal{P}_{\Gamma}\gets \bm{\Gamma} \left(\bm{Y} + \delta\bm{Y}\right) $
            \State $\bm{d}\gets\overline{\mathcal{P}_{\Gamma}}$
            \State Compute $\bm{\Sigma}$ and $\bm{\Delta}$ from $\mathcal{P}_{\Gamma}$
            \State $\bm{\Pi}\gets$BlockCholesky$\left(\bm{\Sigma}\right)$ 
            \State $r\gets +\infty$ 
            \While{$d_{\max}>\varepsilon_{\textrm{N}}\land r>\varepsilon_{\textrm{CV}}$}
                \State $\delta\bm{Y}^*\gets -\bm{\Delta}^{\tr{\textrm{T}}}\tr{\left(\bm{\Pi}^{\tr{\textrm{T}}}\right)^{-1}\bm{\Pi}^{-1} \bm{d}} $
                \State $d_{\max}^*,\ \alpha\gets +\infty,\ 1$ 
                \While{$d_{\max}^*>d_{\max}$}
                    \State $\bm{d}^*\gets\mathcal{P}_{\Gamma}\left(\alpha\delta\bm{Y}^*\right)$
                    \State $d_{\max}^*,\ \alpha\gets \max |\bm{d}^*|, \ \gamma\alpha$
                \EndWhile
                \State $\bm{Y},\ \mathcal{P}_{\Gamma}\gets\bm{Y} + \frac{\alpha}{\gamma}\delta\bm{Y}^*,\ \mathcal{P}_{\Gamma}\left(\frac{\alpha}{\gamma}\delta\bm{Y}^* + \delta\bm{Y}\right)$
                \State $\bm{d}\gets\bm{d}^*$
                \State $r,\ d_{\max}\gets \log d_{\max}^*/\log d_{\max},\ \max |\bm{d}|$
            \EndWhile
        \EndWhile
        \State \textbf{Return:}  $\bm{X}$, $\bm{U}$, $d_{\max}$
    \end{algorithmic}
\end{algorithm}
The acceleration can be estimated in terms of complexity analysis. The complexity for the Cholesky factorisation of a matrix of size $N_c$ is $\frac{N_c^3}{3}$. The complexity of the block Cholesky algorithm for a matrix of $N$ blocks of average size $n_c=\frac{N_c}{N}$ is $\frac{7N_c n_c^2}{3}=\frac{7N_c^3}{3N^2}$. Then, the ratio is $\frac{7}{N^2}$.
Regarding the complexity of the solving of a linear system with backward then forward substitution, the complexity is $N_c^2$ and $n_c^2 + 3(N-1)n_c^2\approx \frac{3N_c^2}{N}$ for the block Cholesky factorisation. Thus, the ratio is $\frac{3}{N}$.
Both the complexity ratios are smaller than \num{1} since the usual number of steps $N$ for astrodynamics problems is generally higher than \num{3}.

\begin{figure}
    \centering
    \tikzstyle{io} = [trapezium, 
trapezium stretches=true, 
trapezium left angle=70, 
trapezium right angle=110, 
thick,
minimum width=3cm, 
minimum height=1cm, text centered, 
draw=phdpurple!50!black, fill=phdpurple!40]

\tikzstyle{process} = [rectangle, 
minimum width=3cm, 
minimum height=1cm, 
text centered, 
text width=3cm, 
thick,
draw=phdblue!50!black, 
fill=phdblue!40]

\tikzstyle{decision} = [diamond, aspect=2.2, 
minimum width=3cm, 
minimum height=1cm, 
text centered,
thick,
draw=phdred2!50!black, 
fill=phdred2!40]
\tikzstyle{arrow} = [thick,->,>=stealth]

\begin{tikzpicture}[node distance=2cm]

\node (inFG) [io] {First guess: $\bm{U}_0$};

\node (proFP) [process, yshift=-0.2cm, below of=inFG] {\textbf{Forward pass}};
\node (decDDP) [decision, below of=proFP] {\makecell[c]{$|J-J^*|\leq\varepsilon_{\textrm{DDP}}$?}};
\node (proBS) [process, right of=decDDP, xshift=2.5cm] {\textbf{Backward sweep}};
\node (DDPalgo) [draw=black,very thick,dotted,inner sep=15pt,
fit=(proFP) (proBS) (proBS)] {};
\node (DDPalgotext) [above of=DDPalgo, xshift=2.6cm, yshift=-0.45cm, text=black] {{\Large\textbf{DDP solver}}};

\node (decAUL) [decision, below of=decDDP, yshift=-0.5cm] {\makecell[c]{$g_{\max}\leq\varepsilon_{\textrm{AUL}}$?}};
\node (proUPAUL) [process, left of=decAUL, xshift=-2.5cm] {
\makecell[c]{Update $\bm{\Lambda}$ and $\bm{M}$}};
\node (proAULcost) [process, left of=proFP, xshift=-2.5cm] {\makecell[c]{Update $\Tilde{l}$ and $\Tilde{\phi}$}};
\node (AULalgo) [draw=black,very thick,dashed,inner sep=4pt, fit=(DDPalgo) (decAUL) (proUPAUL) (proAULcost),align=right,text height=160pt, text=black] {};
\node (AULalgotext) [below of=DDPalgotext, yshift=-3cm, text=black] {{\Large\textbf{AUL solver}}};

\node (decFO) [decision, below of=decAUL, yshift=-0.5cm] {\makecell[c]{$\eta\approx 0$ and $\sigma\approx 0$ ?}};
\node (proUFO) [process, right of=decFO, xshift=2.5cm] {Update $\eta,\sigma$};

\node (proNewton) [process, below of=decFO, yshift=-0.5cm] {\textbf{Newton solver}};

\node (DADDyalgo) [draw=black,ultra thick, fit=(proFP) (decDDP) (proBS) (DDPalgo) (decAUL) (proUPAUL) (proAULcost) (proNewton) (proUFO), align=right, inner sep=10pt] {};
\node (DADDyalgotext) [below of=AULalgotext, yshift=-3cm, text=black] {{\Large\textbf{DADDy solver}}};

\node (out1) [io, below of=proNewton] {Output: $\bm{X}^*, \bm{U}^*, J^*, g_{\max}$};

\draw [arrow] (inFG) -- (proFP);

\draw [arrow] (proFP) -- (decDDP);
\draw [arrow] (decDDP) -- node[anchor=east] {\textbf{Y}} (decAUL);
\draw [arrow] (decDDP) -- node[anchor=south] {\textbf{N}} (proBS);

\draw [arrow] (decAUL) -- node[anchor=east] {\textbf{Y}} (decFO);
\draw [arrow] (decAUL) -- node[anchor=south] {\textbf{N}} (proUPAUL);
\draw [arrow] (proUPAUL) -- (proAULcost);
\draw [arrow] (proAULcost) -- (proFP);

\draw [arrow] (decFO) -- node[anchor=east] {\textbf{Y}} (proNewton);
\draw [arrow] (decFO) -- node[anchor=south] {\textbf{N}} (proUFO);
\draw [arrow] (proUFO.east) -- ++(0.3,0) |- (proFP);

\draw [arrow] (proBS) |- (proFP);
\draw [arrow] (proNewton) -- (out1);

\end{tikzpicture}
    \caption{Summary flow chart of the \gls*{DADDy} solver.}
    \label{fig:flow_chart}
\end{figure}
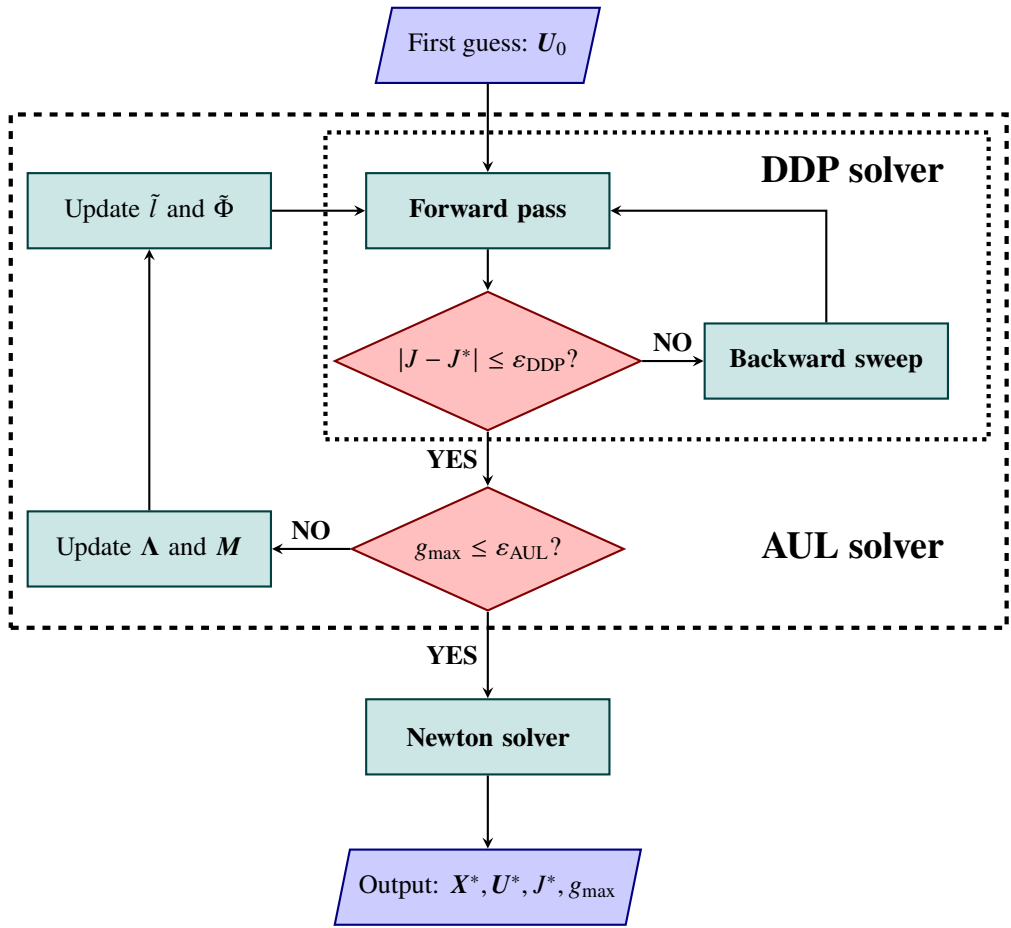

\Fig{fig:flow_chart} represents the implementation of the \gls*{DADDy} solver. This work proposes modifications to the forward pass, the backward sweep, and the Newton solver, which are highlighted in bold in the flow chart.
Two forward pass algorithms and three backward sweep algorithms were presented in \Sec{sec:DA_DDP}. Thus, six different solvers can be tested. They are listed and named in \Tab{tab:DDP_naming} \tr{the \gls*{iLQR} or \gls*{DDP} labels for \Alg{alg:backward_sweep_AD} indicate if the second-order derivatives of the dynamics are taken into account or not}.
\begin{table}[]
    \centering
    \caption{List of the \tr{optimisation} methods and their corresponding acronyms.}
    \begin{tabular}{ c  c  c  c }
        \toprule
         & & \multicolumn{2}{c}{\textbf{Forward pass}} \\ \hline\hline
        & Algorithm & \Alg{alg:forward_pass_AD} & \Alg{alg:forward_pass_DA}  \\ 
        \hline
        \multirow{3}{*}{\textbf{Backward sweep}} & \Alg{alg:backward_sweep_AD} \tr{(\gls*{iLQR})} & \tr{\gls*{iLQR}} & \tr{\gls*{iLQR}}Dyn \\ 
            & \tr{\Alg{alg:backward_sweep_AD} (\gls*{DDP})} & \tr{\gls*{DDP}} & \tr{\gls*{DDP}}Dyn \\ 
            & \Alg{alg:backward_sweep_Q} & Q & QDyn \\
        \bottomrule
    \end{tabular}
    \label{tab:DDP_naming}
\end{table}
\tr{\Gls*{iLQR} and \gls*{DDP} methods respectively corresponds to an implementation of the ALTRO solver \cite{HowellEtAl_2019_AaFSfCTO} and standard \gls*{DDP}} and serve as \tr{standards} for \tr{subsequent} comparisons.
The Newton solver \tr{used} in \gls*{DADDy} is the one developed in \Alg{alg:pn_method_DA}.
\section{Applications} \label{sec:application}
This section focuses on the validation of the \gls*{DADDy} solver and its application to various test cases drawn in the literature. \tr{The following problems are considered:
\begin{enumerate}[label=\Alph*)]
    \item \textbf{Validation \tp{and parameter tuning test case}}: \tp{an Earth–Mars fuel-optimal low-thrust transfer~\citep{LantoineRussell_2012_AHDDPAfCOCPP2A}.}
    \item \textbf{Fuel-optimal low-thrust transfers in the \gls*{CR3BP}}:
    \begin{itemize}
        \item $L_2$ halo to $L_1$ halo~\citep{AzizEtAl_2019_HDDPitCRTBP,BooneMcMahon_2025_RSTOUSTTaDDP}.
        \item $L_2$ \gls*{NRHO} to \gls*{DRO}~\citep{BooneMcMahon_2025_RSTOUSTTaDDP}.
        \item \Gls*{DRO} to \gls*{DRO}~\citep{AzizEtAl_2019_HDDPitCRTBP}.
    \end{itemize}
    \item \tp{\textbf{Fuel-optimal low-thrust transfer in the Geocentric two-body problem}:
    \begin{itemize}
        \item \Gls*{LEO} to \gls*{LEO}~\citep{DiCarloVasile_2021_ASfLTOT}.
        \item \Gls*{MEO} to \gls*{MEO}~\citep{DiCarloVasile_2021_ASfLTOT}.
    \end{itemize}}
\end{enumerate}}
The solver was entirely developed in C++\footnote{Library available at: \url{https://github.com/ThomasClb/DADDy.git} [last accessed \lastdate].}, and uses the \gls*{DACE}\footnote{Library available at: \url{https://github.com/dacelib/dace} [last accessed \lastdate].} \tr{as} polynomial computational engine, implemented by Dinamica SRL for ESA \cite{RasottoEtAl_2016_DASTfNUPiSD,MassariEtAl_2018_DASLwACGfSEA}. All computations and run time analyses were performed on an Intel® Xeon® Gold 6126 CPU at 2.6 GHz. \tp{All settings, except those subjected to the sensitivity analysis ($\varepsilon_{\rm DA}$, $\varepsilon_{\rm AUL}$, and the polynomial order), were kept identical throughout this work to avoid hyper-tuning and are documented in the publicly available code.}

\subsection{Validation \tp{and parameter tuning}.}

\tp{First, we validate the solver and investigate the effects of the various tuning parameters on algorithm performance. The impact of the order of Taylor expansions, the tolerance values, and the algorithm variations are investigated. For all studies performed in this section, we consider the low-thrust Earth-Mars transfer optimisation problem from \citet{LantoineRussell_2012_AHDDPAfCOCPP2A} is considered:
\begin{equation}
    \label{eq:energy_optimal_problem}
    \begin{aligned}
        \bm{x}_k &=\left[\bm{r}_k^{\textrm{T}}, \bm{v}_k^{\textrm{T}}, m_k\right]^{\textrm{T}}, \ \forall k \in \mathbb{N}_{N}, \\
        \bm{f}\left(\bm{x}_k, \bm{u}_k\right) &= \bm{x}_{k} + \int_{t_k}^{t_{k+1}}{\left[\bm{v}^{\tr{\textrm{T}}}, \dot{\bm{v}}^{\tr{\textrm{T}}}, \dot{m}\right]^{\tr{\textrm{T}}}dt}, \ \forall k \in \mathbb{N}_{N-1},\\
        \dot{\bm{v}} &= -\frac{\mu}{\|\bm{r}\|_2^3}\bm{r} + \frac{\bm{u}}{m},\\
        \dot{m} &= -\frac{\|\bm{u}\|_2}{g_0  \tr{\textrm{Isp}}},\\
        \tr{\phi}(\bm{x},\bm{x}_t) &= \left(\bm{r} - \bm{r}_t\right)^{\tr{\textrm{T}}}\left(\bm{r} - \bm{r}_t\right) + \left(\bm{v} - \bm{v}_t\right)^{\tr{\textrm{T}}}\left(\bm{v} - \bm{v}_t\right),\\
        \bm{g}_{\tr{\textrm{ineq}}}(\bm{x},\bm{u}) &= \left[\bm{u}^{\tr{\textrm{T}}}\bm{u} - u_{\tr{\textrm{max}}}^2, m_{\tr{\textrm{dry}}} - m \right]^{\tr{\textrm{T}}},\\
        \bm{g}_{\tr{\textrm{teq}}}(\bm{x},\bm{x}_t) &= \left[\left(\bm{r} - \bm{r}_t\right)^{\tr{\textrm{T}}}, \left(\bm{v} - \bm{v}_t\right)^{\tr{\textrm{T}}}\right]^{\tr{\textrm{T}}},\\
    \end{aligned}
\end{equation}
where $N_x=7$, $N_u=3$, the state vector $\bm{x}$ can be written: $\left[x,y,z,\dot{x}, \dot{y}, \dot{z}, m \right]$ or $\left[\bm{r}^{\tr{\textrm{T}}}, \bm{v}^{\tr{\textrm{T}}}, m \right]^{\tr{\textrm{T}}}$, and $m$ is the spacecraft mass.
The stage cost $\tg{\ell}$ is the same as \Eq{eq:fuel_optimal_ctg}, the number of stages is $N=40$, the \gls*{ToF} is \SI{348.79}{\day}, the initial conditions and target are given in \Tab{tab:ic_earth_mars}.
\begin{table}[]
    \centering
    \caption{Earth-Mars two-body problem transfer initial and target states.}
    \begin{tabular}{l c c c c c c c} 
        \toprule
        \textbf{Type} & \textbf{$x$ [\SI{}{\kilo\meter}]} & \textbf{$y$ [\SI{}{\kilo\meter}]} & \textbf{$z$ [\SI{}{\kilo\meter}]} & \textbf{$\dot{x}$ [\SI{}{\kilo\meter}/\SI{}{\second}]} & \textbf{$\dot{y}$ [\SI{}{\kilo\meter}/\SI{}{\second}]} & \textbf{$\dot{z}$ [\SI{}{\kilo\meter}/\SI{}{\second}]} & \textbf{$m$ [\SI{}{\kilo\gram}]} \\ [0.5ex] 
        \midrule
        Departure & \num{-140699693} & \num{-51614428} & \num{980} & \num{9.774596} & \num{-28.07828} & \num{4.337725E-4} & \num{1000}\\
        Target & \num{-172682023} & \num{176959469} & \num{7948912} & \num{-16.427384} & \num{-14.860506} & \num{9.21486e-2} & – \\
        \bottomrule
    \end{tabular}
    \label{tab:ic_earth_mars}
\end{table}
Normalisation units and various dynamics parameters are reported in \Tab{tab:SUN_units}. All mass parameters in this work were obtained from JPL DE431 ephemerides \footnote{Publicly available at: \url{https://naif.jpl.nasa.gov/pub/naif/generic_kernels/pck/gm_de431.tpc} [retrieved on \lastdate].} \cite{FolknerEtAl_2014_TPaLEDaD}.
\begin{table}[]
    \centering
    \caption{Sun-centred normalisation units and dynamics parameters.}
    \begin{tabular}{ l c c}
        \toprule
      \textbf{Parameter} & \textbf{Symbol} & \textbf{Value} \\
      \midrule
        Mass parameter [-] & $\mu$ & \num{1.32712440041e11}  \\ 
        Length {[\SI{}{\kilo\meter}]} & $ \mathrm{LU} $ & \num{149597870.7}   \\ 
        Time {[\SI{}{\second}]} & $ \mathrm{TU} $ & \num{5022642.891}  \\ 
        Velocity {[\SI{}{\kilo\meter}/\SI{}{\second}]} & $ \mathrm{VU} $ & \num{29.78469183} \\
        Standard gravity [\SI{}{\meter}/\SI{}{\second\squared}] & $g_0$ & \num{9.81}  \\ 
        Specific impulse [\SI{}{\second}] & $\tr{\textrm{Isp}}$ & \num{2000}  \\ 
        Spacecraft dry mass [\SI{}{\kilo\gram}] & $m_{\tr{\textrm{dry}}}$ & \num{500}    \\
        Maximum spacecraft thrust [\SI{}{\newton}] & $u_{\tr{\textrm{max}}}$ & \num{0.5}   \\
        \bottomrule
    \end{tabular}
    \label{tab:SUN_units}
\end{table}
In the remainder of this work, the first guess $\bm{U}_0$ is a $N N_u$ vector with all components equal to \SI{e-6}{\newton}, the initial spacecraft mass, including dry mass, is \SI{1000}{\kilo\gram}, and the various tolerance parameters are set to: $\varepsilon_{\textrm{DDP}}=10^{-4}$, and $\varepsilon_{\textrm{N}}=10^{-10}$.
As discussed in \Sec{sec:fuel_optimal}, fuel-optimal optimisation is performed in four stages for $(\eta,\sigma)$: $(1, 10^{-2}) \rightarrow (0.5, 10^{-2}) \rightarrow (10^{-1}, 2\times10^{-3})\rightarrow (10^{-3}, 10^{-3})$. \Fig{fig:fuel_optimal} shows the solution to this problem, which visually matches the results of \citet{LantoineRussell_2012_AHDDPAfCOCPP2A}. \Fig{fig:fuel_trajectory} shows the trajectory from Earth to Mars and \Fig{fig:fuel_thrust} represents the control norm.
\begin{figure}[H]
    \centering
    \begin{subfigure}{.45\textwidth}
        \centering
        \includegraphics[width=1\linewidth]{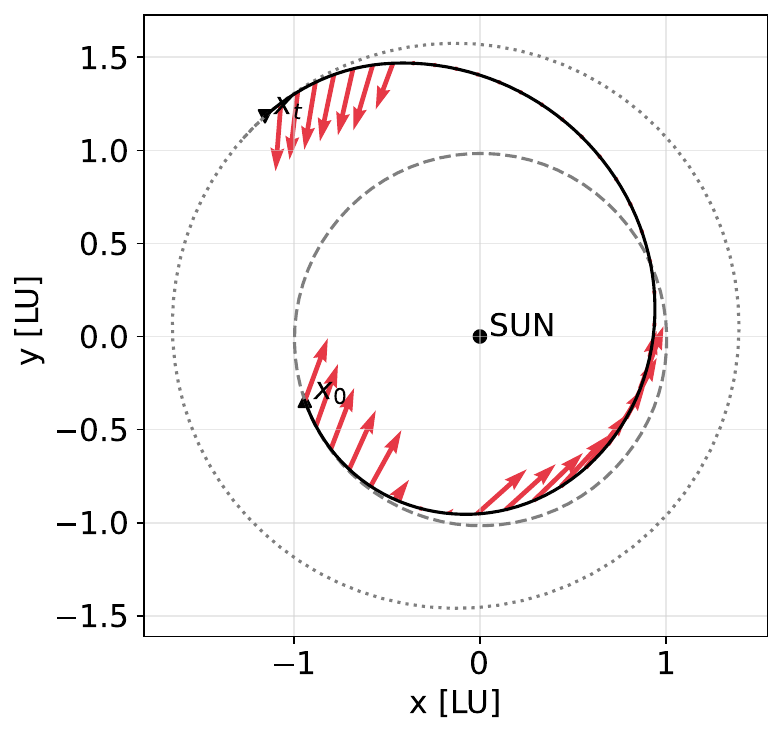}
        \caption{Trajectory in the $x\textrm{--}y$ plane.}
        \label{fig:fuel_trajectory}
    \end{subfigure}%
    \begin{subfigure}{.55\textwidth}
        \centering
        \includegraphics[width=1\linewidth]{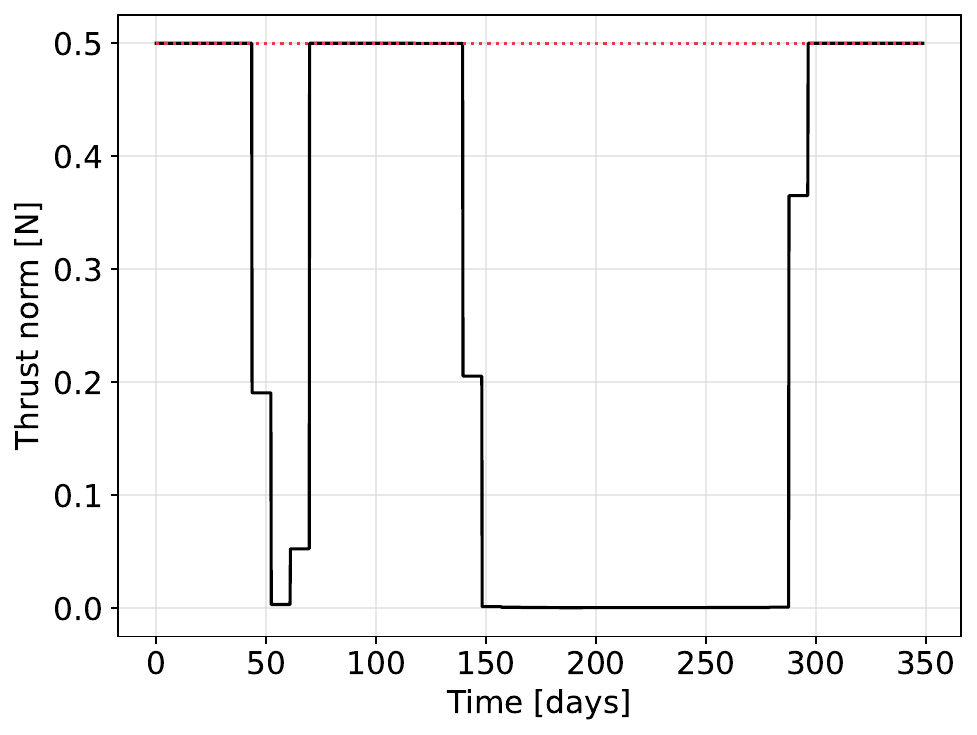}
        \caption{Thrust profile.}
        \label{fig:fuel_thrust}
    \end{subfigure}
    \caption{Solution to the Earth--Mars transfer.}
    \label{fig:fuel_optimal}
\end{figure}}

\subsubsection{\tp{Determining the order of the Taylor expansions}}
\tp{We first investigate the impact of the Taylor polynomial order. Automatic differentiation requires the expansion order to be at least \num{2}. Although higher-order expansions can provide larger convergence radii \citep{ArmellinEtAl_2010_ACECUDAtCoA}, the number of coefficients increases rapidly, thereby significantly raising the complexity of all \gls*{DA} operations \citep{VasileEtAl_2019_SPiDSwGPAaICC}. For this reason, we analyse how the Taylor expansion order affects the run time (RT) and the fuel consumption ($J$). \Tab{tab:res_earth_mars_order} shows the evolution of RT and $J$ for orders \num{2}, \num{3}, and \num{4}, obtained with the iLQRDyn solver using $\varepsilon_{\rm AUL} = 10^{-6} = \varepsilon_{\rm DA}$.}
\begin{table}[h]
    \centering
    \caption{Earth--Mars transfer performance metrics as a function of the expansion order.}
    \begin{tabular}{l c c } 
        \toprule
        \textbf{Order} & \textbf{$J$ [\SI{}{\kilo\gram}]} & \textbf{Run time [\SI{}{\second}]} \\ [0.5ex] 
        \midrule
        2 & \num{396.54} & \num{6.16}\\
        3 & \num{396.54} & \num{18.2} \\
        4 & \num{396.56} & \num{84.1} \\
        \bottomrule
    \end{tabular}
    \label{tab:res_earth_mars_order}
\end{table}

\tp{The results indicate that all configurations converge to the same solution. However, significant differences are observed in run time: lower orders lead to faster convergence. This behaviour is consistent with the findings of \citet{BooneMcMahon_2025_RSTOUSTTaDDP}, where lower orders also reduce run time but at the expense of solution quality. The key distinction is that \citet{BooneMcMahon_2025_RSTOUSTTaDDP} employ a single Taylor expansion throughout the entire re-optimisation process, whereas in the present work the trajectory is automatically re-expanded whenever the previous expansion becomes obsolete. Based on these results, the expansion order is set to the smallest possible value, \ie \num{2}.}

\subsubsection{\tp{Setting the tolerance of the dynamics approximation trigger}}
\tp{We now examine the impact of the tolerance parameter $\varepsilon_{\rm DA}$. This parameter governs the decision of whether the dynamics should be updated using dynamics approximation with Taylor expansions or recomputed from scratch.
\Fig{fig:res_earth_mars_DA_tol} illustrates the evolution of RT and $J$ as $\varepsilon_{\rm DA}$ varies with the iLQRDyn solver using $\varepsilon_{\rm AUL} = 10^{-6}$. Results are normalised with respect to the case $\varepsilon_{\rm DA}=10^{-6}$.}
\begin{figure}[H]
    \centering
    \begin{tikzpicture}
\begin{axis}[
    xlabel={$\varepsilon_{\rm DA}$},
    xmode=log, 
    width=7cm,
    height=7cm,
    ylabel={Value compared to $\varepsilon_{\rm DA}=10^{-6}$ [\%]},
    ymin=0.75, ymax=1.35,
    legend style={at={(0.9,0.9)}, anchor=north east}, 
]

\addplot[phdred2, thick, mark=*] coordinates {
    (1e-4, 0.818)
    (5e-5, 0.899)
    (2e-5, 0.862)
    (1e-5, 0.939)
    (5e-6, 0.966)
    (2e-6, 1.003)
    (1e-6, 1)
    (5e-7, 1.085)
    (2e-7, 1.021)
    (1e-7, 1.171)
    (5e-8, 1.158)
    (2e-8, 1.205)
    (1e-8, 1.301)
};
\addlegendentry{RT}

\addplot[phdblue, thick, mark=square*] coordinates {
    (1e-4, 1.000)
    (5e-5, 1.000)
    (2e-5, 1.000)
    (1e-5, 1.000)
    (5e-6, 1.000)
    (2e-6, 1.000)
    (1e-6, 1)
    (5e-7, 1.000)
    (2e-7, 1.000)
    (1e-7, 1.000)
    (5e-8, 1.000)
    (2e-8, 1.000)
    (1e-8, 1.000)
};
\addlegendentry{$J$}

\end{axis}
\end{tikzpicture}
    \caption{Earth--Mars transfer performance metrics for different values of $\varepsilon_{\rm DA}$.}
    \label{fig:res_earth_mars_DA_tol}
\end{figure}
\tp{The results show that the cost functions remain identical up to four decimal places, while the run time increases as $\varepsilon_{\rm DA}$ decreases. This behaviour is consistent with the fact that larger values of $\varepsilon_{\rm DA}$ lead to more frequent reuse of the dynamics approximation, thereby reducing computational effort. However, using too large a value of $\varepsilon_{\rm DA}$ may result in poor approximations of the dynamics, which in turn can introduce errors and constraint violations.
For consistency, and to ensure that constraint violations are not caused by inaccuracies in the Taylor model, $\varepsilon_{\rm DA}$ is set equal to $\varepsilon_{\rm AUL}$ throughout this work.}

\subsubsection{\tp{Tolerance selection for the Newton‑solver trigger}}
\tp{We now analyse the impact of the tolerance parameter $\varepsilon_{\rm AUL}$, which determines when the \gls*{AUL} solver stops and the Newton solver is triggered.
To assess its influence, we compare the run time, cost function, number of \gls*{iLQR}/\gls*{DDP} iterations ($n$ DDP), number of \gls*{AUL} iterations ($n$ AUL), and number of Newton solver iterations ($n$ Newton) for values of $\varepsilon_{\rm AUL}$ ranging from \num{e-6} to \num{e-10}, corresponding to the tolerance used by the Newton solver. Note that in the latter case, the Newton solver is not triggered, as the constraint violation is already satisfactory upon exiting the AUL solver.
\Tab{tab:res_earth_mars_AUL_tol} illustrates the evolution of these quantities as $\varepsilon_{\rm AUL}$ varies with the iLQRDyn solver.}
\begin{table}[h]
    \centering
    \caption{Earth--Mars transfer performance metrics for different values of $\varepsilon_{\rm AUL}$.}
    \begin{tabular}{l c c c c c} 
        \toprule
        \textbf{$\varepsilon_{\rm AUL}$} & \textbf{$J$ [\SI{}{\kilo\gram}]} & \textbf{Run time [\SI{}{\second}]} & \textbf{$n$ DDP} & \textbf{$n$ AUL} & \textbf{$n$ Newton}\\ [0.5ex] 
        \midrule
        \num{e-2} & \num{398.22} & \num{4.30} & \num{298} & \num{4} & \num{35} \\
        \num{e-4} ($\varepsilon_{\rm DDP}$) & \num{396.53} & \num{4.42} & \num{353} & \num{20} & \num{44}\\
        \num{e-6} & \num{396.54} & \num{6.16} & \num{454} & \num{35} & \num{11}\\
        \num{e-8} & \num{396.55} & \num{9.44} & \num{492} & \num{41} & \num{38}\\
        \num{e-10} ($\varepsilon_{\rm N}$) & \num{396.54} & \num{11.58} & \num{759} & \num{216} & \num{0} \\
        \bottomrule
    \end{tabular}
    \label{tab:res_earth_mars_AUL_tol}
\end{table}
\tp{The results indicate that activating the Newton solver earlier improves performance, as reflected by reduced run times and fewer iterations for both the \gls*{DDP} and \gls*{AUL} solvers. This observation aligns with the conclusions of \citet{Howell_1984_TDPhO}, since the Newton method converges much faster than \gls*{iLQR}/\gls*{DDP}. The cost function remains unchanged up to four significant digits, except for the case $\varepsilon_{\rm AUL} = 10^{-2}$, where the solution consumes \SI{1.68}{\kilo\gram} more fuel than the $\num{e-10}$ reference (from \SI{396.54}{\kilo\gram} to \SI{398.22}{\kilo\gram}), corresponding to an increase of about \SI{0.4}{\percent}.  
This discrepancy arises because the Newton method is neither an optimisation solver nor a global algorithm. If the constraints are already satisfied for $\varepsilon_{\rm AUL} > \varepsilon_{\rm DDP}$, polishing the solution with Newton iterations may alter the cost function and yield a sub-optimal result, even though the constraints are enforced with high precision $\varepsilon_{\rm N} \ll \varepsilon_{\rm DDP}$. Conversely, when $\varepsilon_{\rm AUL} \ll \varepsilon_{\rm DDP}$, the Newton method affects the cost by less than the optimisation tolerance, \ie below the convergence criterion of the solver.  
Moreover, the Newton method is local and highly sensitive to the quality of the initial guess, which motivates choosing $\varepsilon_{\rm AUL}$ sufficiently small to ensure reliable convergence. Based on these considerations, the Newton-solver trigger is set to $\varepsilon_{\rm AUL} = \varepsilon_{\rm DDP}/100 = 10^{-6}$ in this work. This choice guarantees that the Newton solver converges robustly without significantly altering the optimal solution while still reducing the run time by nearly a factor of two. It also allows each solver to operate in its most effective regime: \gls*{iLQR}/\gls*{DDP} finds an optimal, nearly feasible fuel-optimal solution from scratch, and the Newton method rapidly converges to a fully feasible solution.}

\subsubsection{\tp{Performance metrics}}
The convergence results are reported in \Tab{tab:DADDy_res_earth_mars_fuel}.
\begin{table}[]
    \centering
    \caption{Fuel-optimal Earth-Mars transfer performance metrics.}
    \footnotesize
    \begin{tabular}{l c c c c c c} 
        \toprule
        \textbf{Data} & \textbf{ \gls*{iLQR}} & \textbf{ \gls*{DDP}} & \textbf{Q} & \textbf{ \gls*{iLQR}Dyn} & \textbf{ \gls*{DDP}Dyn} &\textbf{QDyn} \\ [0.5ex] 
        \midrule
        $J$ [\SI{}{\kilo\gram}] & \num{396.68} & \num{396.78} & \num{397.22} & \num{396.56} & \num{396.74} & \num{396.62} \\
        Run time [\SI{}{\second}] & \num{27} & \num{18} & \num{16} & \num{8.2} & \num{8.4} & \num{8.5} \\
        \bottomrule
    \end{tabular}
    \label{tab:DADDy_res_earth_mars_fuel}
\end{table}
All methods reach satisfactory constraints violation and converge to the same trajectory with three significant digits. 
The \tr{\gls*{DDP}} and Q methods reduce overall runtime. 
Yet, these gains vanish with dynamics approximation, i.e., for \tr{\gls*{DDP}}Dyn and QDyn.
\tp{Indeed, the three methods that employ dynamics approximations achieve similar run times and reduce the computational burden by between \SI{47}{\percent} and \SI{70}{\percent} compared to their counterparts that do not.}
\tr{This supports the findings of \citet{BooneMcMahon_2025_RSTOUSTTaDDP}, who observed that dynamics evaluation accounts for a significant portion of the total runtime.} \tp{Reducing the computational requirements for these repeated evaluations can lead to significant runtime improvements for the overall algorithm.}

\begin{figure}[H]
    \centering
    \includegraphics[width=0.5\linewidth]{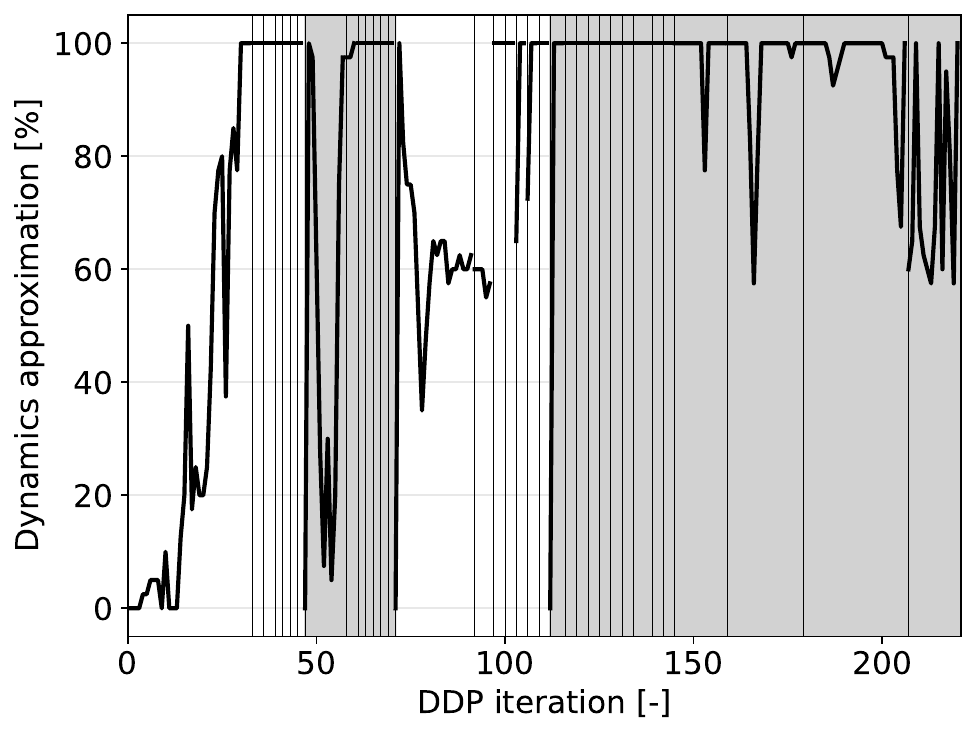}
    \caption{Proportion of dynamics approximations performed by the \tr{\gls*{iLQR}}Dyn method during the Earth--Mars transfer.}
    \label{fig:dynamics_approximation}
\end{figure}

\tr{\Fig{fig:dynamics_approximation} shows the proportion of dynamics propagation performed using polynomial approximation at each iteration for the Earth-Mars fuel-optimal transfer using the method \gls*{iLQR}Dyn. A value of \num{100}\% indicates that all $N$ state propagations were carried out using polynomial expansions, whereas a value of \num{0}\% means that all $N$ states were computed from scratch.
Vertical lines mark \tp{updates of the dual states $\bm{\Lambda}$ and the penalty factors $\bm{M}$, \ie the \gls*{AUL} solver iterations.}
The alternating gray and white regions represent successive stages of the fuel-optimal optimisation process: the first white region corresponds to the energy-optimal optimisation with $(\eta,\sigma) = (1, 10^{-2})$, the first gray region to the phase with $(\eta,\sigma) = (0.5, 10^{-2})$, and so on.}
\tr{This figure shows that the dynamics are primarily computed from scratch at the beginning of the energy-optimal optimisation phase and at the start of each update to the pair $(\eta,\sigma)$, when the trajectory undergoes significant changes.
Overall, during the entire fuel-optimal optimisation process, the dynamics approximation is used on average \num{79.4}\% of the time. This high approximation rate contributes to the substantial reduction in run time observed in \Tab{tab:DADDy_res_earth_mars_fuel} for the methods that implement dynamics approximation during the forward pass: \gls*{iLQR}Dyn, \gls*{DDP}Dyn, and QDyn.}

\subsection{Earth-Moon \gls*{CR3BP} transfers}
The solver is now tested on trajectory optimisation problems in the Earth-Moon \gls*{CR3BP} system \citep{Poincare_1892_LMNDLMC,Szebehely_1967_ToO}. The dynamics take the following form:
\begin{equation}
    \begin{aligned}
        \bm{f}(\bm{x},\bm{u}) & =  \bm{x}_{k} + \int_{t_k}^{t_{k+1}}{\left[\dot{x}, \dot{y}, \dot{z},\Ddot{x} ,\Ddot{y} ,\Ddot{z}, \dot{m}\right]^{\tr{\textrm{T}}}}dt,\\
        \Ddot{x} & =  2\dot{y} + \dfrac{\partial\Omega}{\partial x} + \dfrac{u^x}{m},\\
        \Ddot{y} &  = -2\dot{x} + \dfrac{\partial\Omega}{\partial y} + \dfrac{u^y}{m},\\
        \Ddot{z} & =  \dfrac{\partial\Omega}{\partial z} + \dfrac{u^z}{m},\\
        \dot{m} &= -\frac{\|\bm{u}\|_2}{g_0 \tr{\textrm{Isp}}},\\
    \end{aligned}
\end{equation}
where $N_x=7$, $N_u=3$, $\bm{x}=\left[x,y,z,\dot{x}, \dot{y}, \dot{z}, m \right]=\left[\bm{r}^{\tr{\textrm{T}}}, \bm{v}^{\tr{\textrm{T}}}, m \right]^{\tr{\textrm{T}}}$, $\bm{u}=\left[u^x,u^y,u^z\right]$, and $\Omega =  \frac{1}{2}\left(x^2 + y^2\right) + \dfrac{1 - \mu}{\sqrt{(x+\mu)^2 + y^2 + z^2}} + \dfrac{\mu}{\sqrt{(x+\mu-1)^2 + y^2 + z^2}}$. The normalisation units \cite{JorbaMasdemont_1999_DitCMotCPotRTBP} are reported in \Tab{tab:EARTH_MOON_units}.
\begin{table}[]
    \centering
    \caption{Earth-Moon \gls*{CR3BP} normalisation units and parameters.}
    \begin{tabular}{ l l l}
        \toprule
      \textbf{Parameter} & \textbf{Symbol} & \textbf{Value} \\
      \midrule
        Mass parameter of $M_1$ [\SI{}{\kilo\meter\cubed}/\SI{}{\second\squared}] & $ \mathrm{GM}_1 $ & \SI{398600}{} \\
        Mass parameter of $M_2$ [\SI{}{\kilo\meter\cubed}/\SI{}{\second\squared}] & $ \mathrm{GM}_2 $ & \SI{4902.80}{}  \\
        Mass parameter {[-]}& $ \mu $ & \SI{1.21506e-2}{}  \\
        Length [\SI{}{\kilo\meter}] & $ \mathrm{LU} $ & \SI{384399}{}  \\ 
        Time [\SI{}{\second}] & $ \mathrm{TU} $ & \SI{375189}{} \\
        Velocity [\SI{}{\kilo\meter}/\SI{}{\second}] & $ \mathrm{VU} $ & \SI{1.02455}{}  \\
        \bottomrule
    \end{tabular}
    \label{tab:EARTH_MOON_units}
\end{table}
The stage cost is the one from \Eq{eq:fuel_optimal_ctg}, and the terminal cost, the path constraints and the terminal constraints are the same as in \Eq{eq:energy_optimal_problem}. Moreover, the spacecraft parameters are similar to those of \Tab{tab:SUN_units} and the tolerances are the same.
The solver was tested on three \gls*{CR3BP} test cases:
\tr{\begin{enumerate}
    \item A transfer from a $L_2$ halo orbit \cite{Farquhar_1970_TCaUoLPS,Howell_1984_TDPhO} to a $L_1$ halo \tp{inspired} by \citet{AzizEtAl_2019_HDDPitCRTBP} and \citet{BooneMcMahon_2025_RSTOUSTTaDDP}.
    \item A transfer from a $L_2$ \gls*{NRHO} \citep{ZimovanSpreenEtAl_2020_NRHOaNHPDSOSaRP} to a \gls*{DRO} \cite{Henon_1969_NEotRPV} inspired by \citet{BooneMcMahon_2025_RSTOUSTTaDDP}.
    \item A \gls*{DRO} to \gls*{DRO} transfer from \citet{AzizEtAl_2019_HDDPitCRTBP}.
\end{enumerate}}
The initial conditions, targets, \gls*{ToF}s and number of stages for each transfer are given in \Tab{tab:ic_cr3bp}.
The values for the $\tr{\textrm{Isp}}$, $g_0$, the maximum thrust magnitude, and the dry mass of the spacecraft are the same as in the Earth-Mars transfer test case, reported in \Tab{tab:SUN_units}. \tp{The value of the maximum thrust magnitude is divided by \num{5} for the \gls*{DRO} to \gls*{DRO} transfer, similarly to \citet{AzizEtAl_2019_HDDPitCRTBP}.}
\begin{table}[]
    \centering
    \caption{Earth-Moon \gls*{CR3BP} transfers data.}
    \begin{tabular}{l c c c c c c c c c} 
        \toprule
        \textbf{Transfer} & \textbf{\gls*{ToF}} [\SI{}{\day}] & $N$ & \textbf{Type} & \textbf{$x$ [LU]} & \textbf{$y$ [LU]} & \textbf{$z$ [LU]} & \textbf{$\dot{x}$ [VU]} & \textbf{$\dot{y}$ [VU]} & \textbf{$\dot{z}$ [VU]} \\ [0.5ex] 
        \midrule
        \multirow{2}{*}{\parbox{2.2cm}{Halo $L_2$ to halo $L_1$}} & \multirow{2}{*}{\tp{\num{32.25}}} & \multirow{2}{*}{\tp{\num{150}}} & Initial & \num{1.16080} & \num{0} & \num{-0.12270} & \num{0} & \num{-0.20768} & \num{0} \\
        & & & Target & \num{0.84871} & \num{0} & \num{0.17389} & \num{0} & \num{0.26350} & \num{0} \\
        \midrule[0.5pt]
        \multirow{2}{*}{\parbox{2.2cm}{\Gls*{NRHO} $L_2$ to \gls*{DRO}}} & \multirow{2}{*}{\num{21.2}} & \multirow{2}{*}{\num{150}} & Initial & \num{1.02197} & \num{0} & \num{-0.18206} & \num{0} & \num{-0.10314} & \num{0} \\
        & & & Target & \num{0.98337} & \num{0.25921} & \num{0} & \num{0.35134} & \num{-0.00833} & \num{0} \\
        \midrule[0.5pt]
        \multirow{2}{*}{\parbox{2.0cm}{\gls*{DRO} to \gls*{DRO}}} & \multirow{2}{*}{\tp{\num{51.25}}} & \multirow{2}{*}{\tp{\num{100}}} & Initial & \num{1.17136} & \num{0} & \num{0} & \num{0} & \num{-0.48946} & \num{0} \\
        & & & Target & \num{1.30184} & \num{0} & \num{0} & \num{0} & \num{-0.64218} & \num{0} \\
        \bottomrule
    \end{tabular}
    \label{tab:ic_cr3bp}
\end{table}

\begin{figure}[H]
    \centering 
    \begin{subfigure}{0.29\textwidth}
      \centering \includegraphics[width=\linewidth]{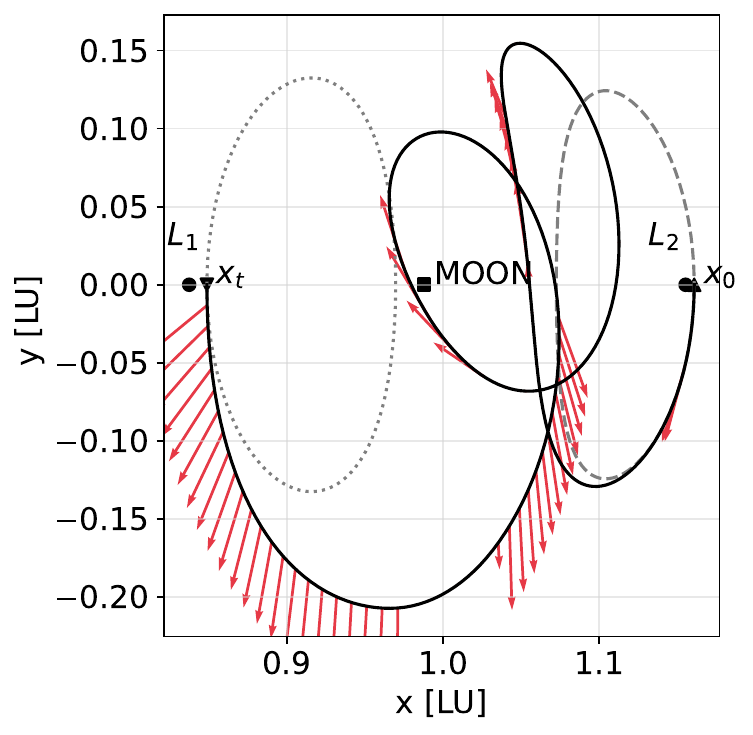}
      \caption{Trajectory in the $x\textbf{--}y$ plane.}
      \label{fig:halo_to_halo_xy}
    \end{subfigure}\hfil 
    \begin{subfigure}{0.28\textwidth}
      \centering \includegraphics[width=\linewidth]{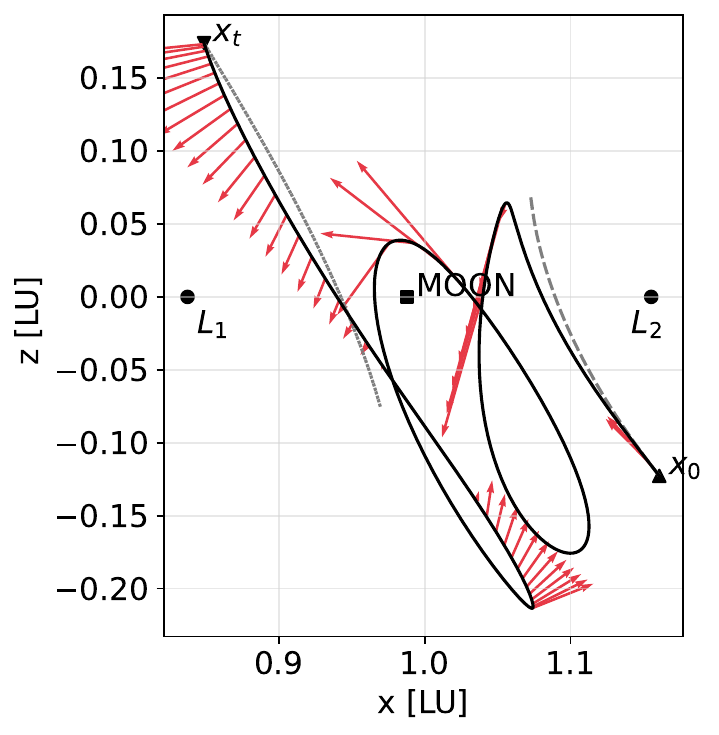}
      \caption{Trajectory in the $x\textbf{--}z$ plane.}
      \label{fig:halo_to_halo_xz}
    \end{subfigure}
    \hfil 
    \begin{subfigure}{0.38\textwidth}
      \centering 
      \includegraphics[width=\linewidth]{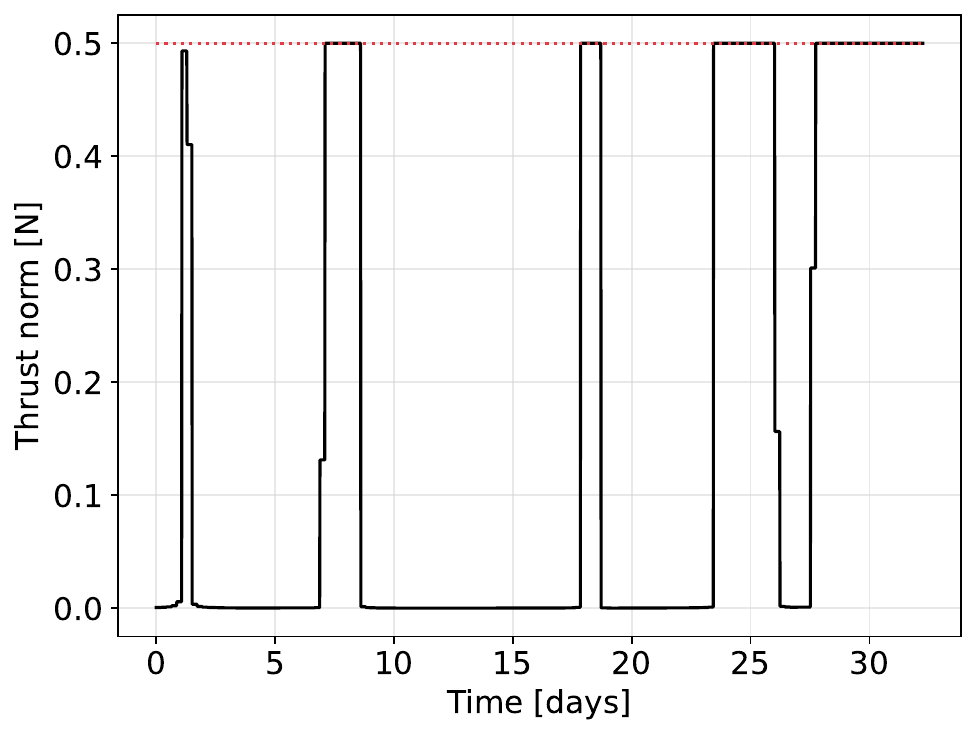}
      \caption{Thrust profile.}
      \label{fig:halo_to_halo_thrust}
    \end{subfigure}
    \caption{Solution to the halo $L_2$ to halo $L_1$ transfer.}
    \label{fig:halo_to_halo}
\end{figure}

\begin{figure}[H]
    \centering
    \begin{subfigure}{.35\textwidth}
        \centering
        \includegraphics[width=1\linewidth]{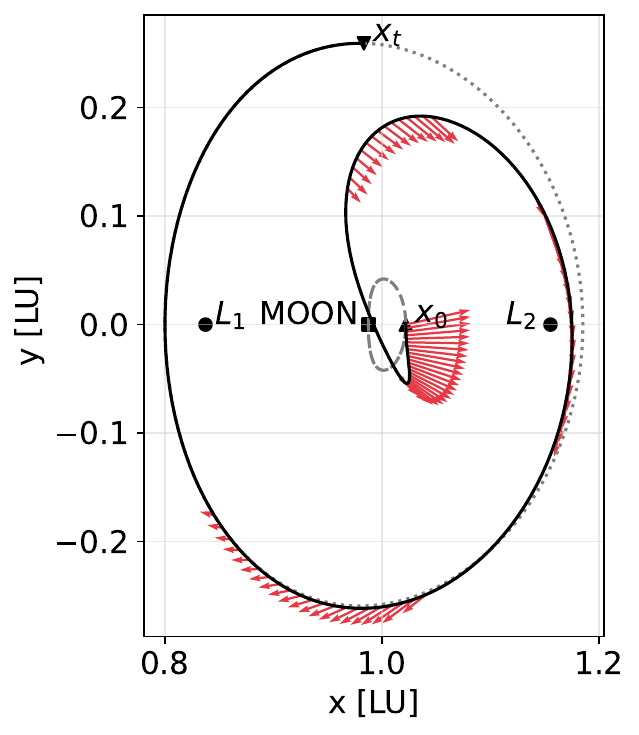}
        \caption{Trajectory in the $x\textbf{--}y$ plane.}
        \label{fig:nrho_to_dro_trajectory}
    \end{subfigure}%
    \begin{subfigure}{.54\textwidth}
        \centering
        \includegraphics[width=\linewidth]{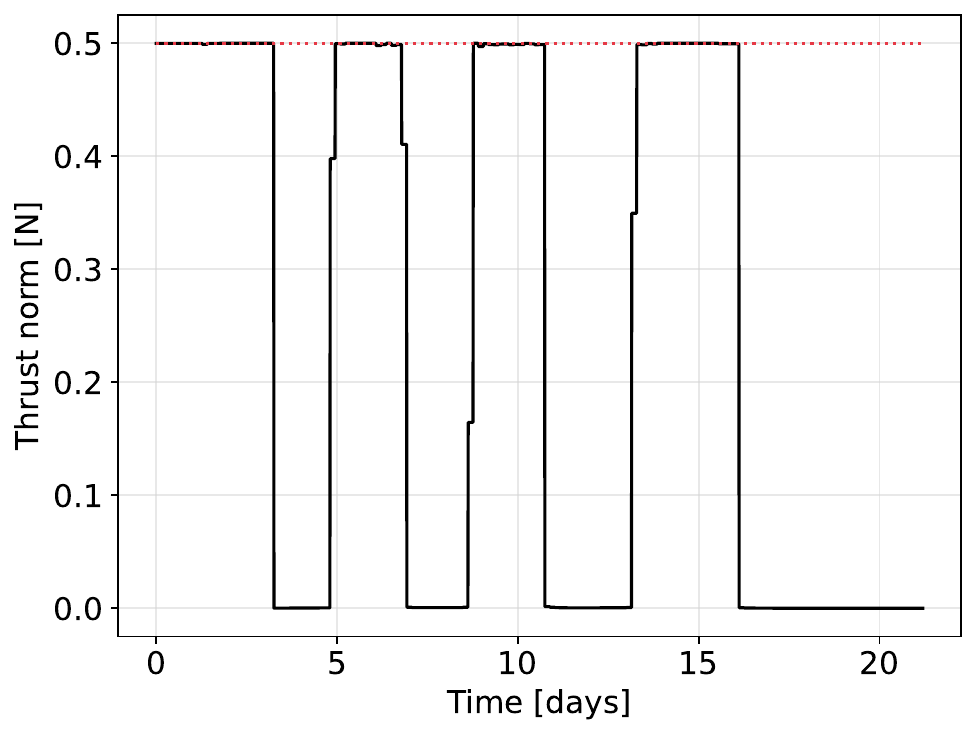}
        \caption{Thrust profile.}
        \label{fig:nrho_to_dro_thrust}
    \end{subfigure}
    \caption{Solution to the $L_2$ \gls*{NRHO} to \gls*{DRO} transfer.}
    \label{fig:nrho_to_dro}
\end{figure}

\begin{figure}[H]
    \centering 
    \begin{subfigure}{0.37\textwidth}
      \centering 
      \includegraphics[width=0.85\linewidth]{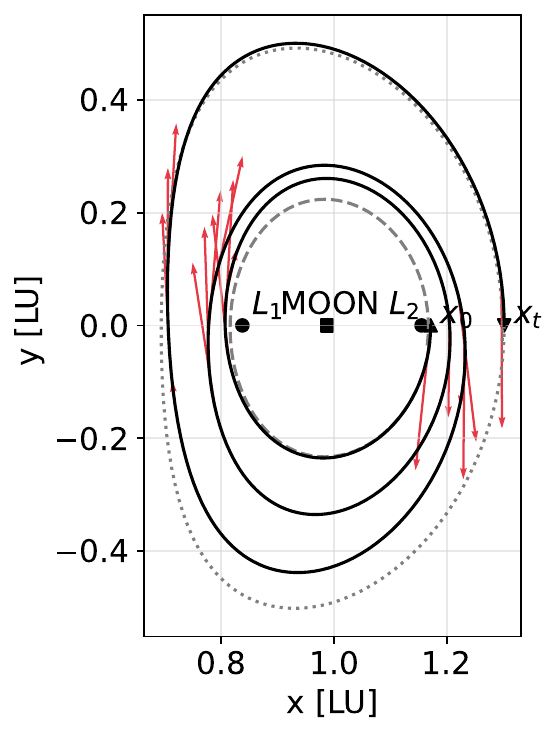}
      \caption{Trajectory in the $x\textbf{--}y$ plane.}
      \label{fig:dro_to_dro_trajectory}
    \end{subfigure}\hfil 
    \begin{subfigure}{0.63\textwidth}
      \centering 
      \includegraphics[width=0.85\linewidth]{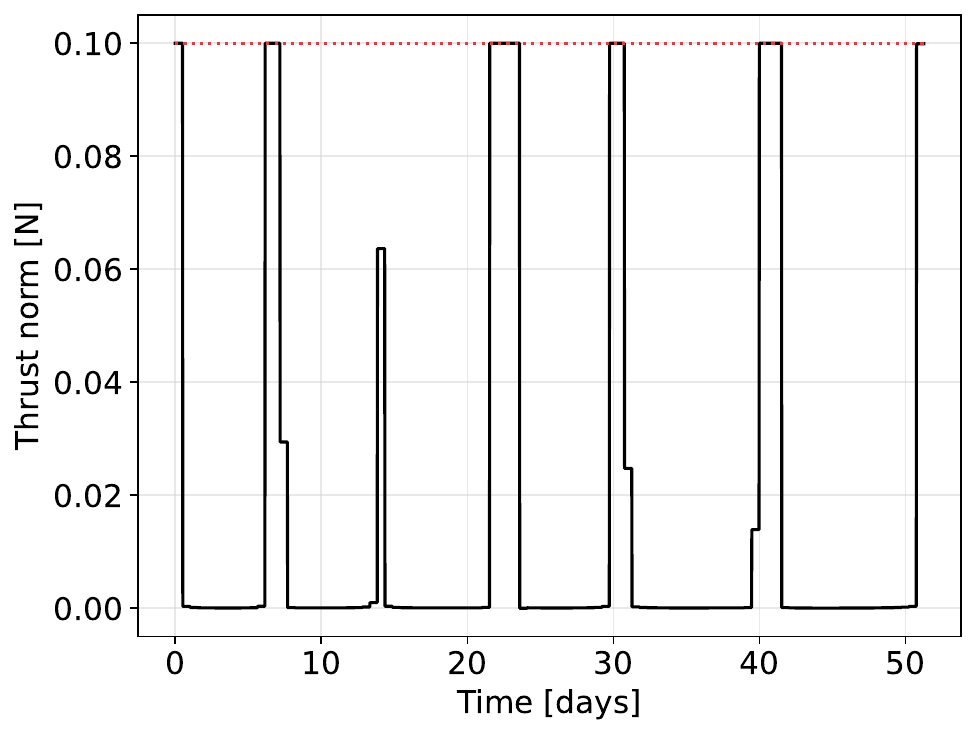}
      \caption{Thrust profile.}
      \label{fig:dro_to_dro_thrust}
    \end{subfigure}
    \caption{Solution to the \gls*{DRO} to \gls*{DRO} transfer.}
    \label{fig:dro_to_dro}
\end{figure}

\Fig{fig:halo_to_halo} shows the solutions to the halo $L_2$ to halo $L_1$ fuel-optimal transfers. \Fig{fig:halo_to_halo_xy}, \Fig{fig:halo_to_halo_xz}, and \Fig{fig:halo_to_halo_thrust} show respectively the transfer in the $x\textbf{--}y$ plane, in the $x\textbf{--}z$ plane and its thrust profile, which is similar to \citet{AzizEtAl_2019_HDDPitCRTBP} and \citet{BooneMcMahon_2025_RSTOUSTTaDDP}. 
\Fig{fig:nrho_to_dro} shows the solution to the \gls*{NRHO} to \gls*{DRO} fuel-optimal transfer. \Fig{fig:nrho_to_dro_trajectory} shows the trajectory in the $x\textbf{--}y$ plane, while \Fig{fig:nrho_to_dro_thrust} shows the thrust profile.
\tg{Finally, \Fig{fig:dro_to_dro} shows the solution to the \gls*{DRO} to \gls*{DRO} fuel-optimal transfer, which correspond to the results of \citet{AzizEtAl_2019_HDDPitCRTBP}. 
\Fig{fig:dro_to_dro_trajectory} and \Fig{fig:dro_to_dro_thrust} respectively show the trajectory and the thrust profile.}

As expected \tr{for fuel-optimal transfers}, all resulting solutions correspond to bang-bang control laws. The performance metrics are given in \Tab{tab:DADDy_res_cr3bp}, where the mention \gls*{DNC} is given when the algorithm did not terminate.
\begin{table}[H]
    \centering
    \caption{Performance metrics for fuel-optimal \gls*{CR3BP} transfers.}
    \scriptsize
    \begin{tabular}{l l c c c c c c} 
        \toprule
        \textbf{Transfer} & \textbf{Data} & \textbf{ \gls*{iLQR}} & \textbf{ \gls*{DDP}} & \textbf{Q} & \textbf{ \gls*{iLQR}Dyn} & \textbf{ \gls*{DDP}Dyn} &\textbf{QDyn} \\ [0.5ex] 
        \midrule
        \multirow{2}{*}{\parbox{1.8cm}{Halo $L_2$ to halo $L_1$}} 
        & $J$ [\SI{}{\kilo\gram}] & \num{26.69} & \num{26.17} & \num{26.22} & \num{26.07} & \num{26.20} & \num{26.23}\\
        & Run time [\SI{}{\second}] & \num{110}  & \num{57}  & \num{57} & \num{31}  & \num{28} & \num{30} \\
        \midrule[0.5pt]

        \multirow{2}{*}{{\parbox{1.8cm}{\Gls*{NRHO} $L_2$ to \gls*{DRO}}}}
        & $J$ [\SI{}{\kilo\gram}] & \num{22.66} & \multicolumn{2}{c}{\multirow{2}{*}{\Gls*{DNC}}} & \num{22.61} & \multicolumn{2}{c}{\multirow{2}{*}{\Gls*{DNC}}} \\
        & Run time [\SI{}{\second}] & \num{410} &  &  & \num{50}  &  & \\
        \midrule[0.5pt]

        \multirow{2}{*}{{\parbox{1.8cm}{\gls*{DRO} to \gls*{DRO}}}}
        & $J$ [\SI{}{\kilo\gram}] & \num{3.26} & \num{5.16} & \num{5.20} & \num{3.25} & \num{5.15} & \multirow{2}{*}{\Gls*{DNC}}\\
        & Run time [\SI{}{\second}] & \num{112}& \num{161} & \num{236} & \num{46} & \num{86} & \\
        \bottomrule
    \end{tabular}
    \label{tab:DADDy_res_cr3bp}
\end{table}
\tp{Results show that DDP, DDPDyn, Q, QDyn, and Dyn versions converge in less than half of the test cases.
This observation aligns with the findings of \citet{NgangaWensing_2021_ASODDPfRBS}, who note that while \gls*{DDP} can converge more rapidly than \gls*{iLQR}, it often requires additional regularisation to ensure stable convergence. Conversely, the \gls*{iLQR}Dyn method consistently delivers results comparable to \gls*{iLQR} with 51\%–88\% shorter run times. This trend holds for all “Dyn” methods, which uniformly outperform their classical counterparts by a substantial margin in computational efficiency while achieving similar cost‑function values. Additionally, although the underlying cause remains unclear, the Q and QDyn methods are never more efficient or stable than the \gls*{DDP} and \gls*{DDP}Dyn methods.}

\subsection{Earth-centred transfer}
Earth-centred two-body \cite{ValladoMcClain_2007_FoAaA} optimisation problems were also solved.
The Gauss equations of motion are written in the equinoctial form from \citet{DiCarloEtAl_2021_EAFftPKMuLTAaOP}:
\begin{equation}
    \begin{aligned}
        \bm{f}\left(\bm{x}_k, \bm{u}_k\right) &= \bm{x}_{k} + \int_{t_k}^{t_{k+1}}{\left[\dot{a}, \dot{p}, \dot{q},\dot{r} ,\dot{s} ,\dot{\tr{L}}, \dot{m}\right]^{\tr{\textrm{T}}}}dt,\\
        \dot{a} & = \frac{2}{\mathcal{B}}\sqrt{\frac{a^3}{\mu}}\left[\left(q\sin \tr{L} - p \cos \tr{L}\right) \frac{u^{\rm R}}{\tp{m}}+\Psi\right],\\
        \dot{p} & = \mathcal{B}\sqrt{\frac{a}{\mu}}\left[-\cos \tr{L}  \frac{u^{\rm R}}{\tp{m}} + \left(\frac{p + \sin \tr{L}}{\Psi} +  \sin \tr{L}\right) \frac{u^{\rm T}}{\tp{m}} -  q\frac{r\cos \tr{L} - s\sin \tr{L}}{\Psi} \frac{u^{\rm N}}{\tp{m}} \right],\\
        \dot{q} & = \mathcal{B}\sqrt{\frac{a}{\mu}}\left[\sin \tr{L}  \frac{u^{\rm R}}{\tp{m}} + \left(\frac{q + \cos \tr{L}}{\Psi} +  \cos \tr{L}\right) \frac{u^{\rm T}}{\tp{m}} +  p\frac{r\cos \tr{L} - s\sin \tr{L}}{\Psi} \frac{u^{\rm N}}{\tp{m}} \right],\\
        \dot{r} & = \frac{\mathcal{B}}{2}\sqrt{\frac{a}{\mu}}\left(1+r^2 + s^2 \right) \frac{\sin \tr{L}}{\Psi} \frac{u^{\rm N}}{\tp{m}},\\
        \dot{s} & = \frac{\mathcal{B}}{2}\sqrt{\frac{a}{\mu}}\left(1+r^2 + s^2 \right) \frac{\cos \tr{L}}{\Psi} \frac{u^{\rm N}}{\tp{m}},\\
        \dot{\tr{L}} & = \sqrt{\frac{\mu}{a}}\frac{\Psi^2 }{ \mathcal{B}^3} - \sqrt{\frac{a^3}{\mu}}\frac{\mathcal{B} }{\Psi}\left(r\cos \tr{L} - s \sin \tr{L}\right) \frac{u^{\rm N}}{\tp{m}} ,\\
        \dot{m} &= -\frac{\|\bm{u}\|_2}{g_0 \tr{\textrm{Isp}}},\\
        \tr{\phi}(\bm{x},\bm{x}_t) &= \left(\bm{r} - \bm{r}_t\right)^{\tr{\textrm{T}}}\left(\bm{r} - \bm{r}_t\right),\\
        \bm{g}_{\tr{\textrm{teq}}}(\bm{x},\bm{x}_t) &= \bm{r} - \bm{r}_t,\\
    \end{aligned}
\end{equation}
where $N_x=7$, $N_u=3$, $\mathcal{B}=\sqrt{1-p^2-q^2}$, and $\Psi=1+p\sin \tr{L} + q\cos \tr{L}$.
\tr{The equinoctial coordinates $\left[a,p,q,r,s,\tr{L} \right]^{\tr{\textrm{T}}}=\left[\bm{r}^{\tr{\textrm{T}}}, \tr{L} \right]^{\tr{\textrm{T}}}$ are defined as:} 
\begin{equation}
    \begin{aligned}
        a, & \\
        p &= e \sin \left(\Omega + \omega\right), \\
        q &= e \cos \left(\Omega + \omega\right), \\
        r &= \tan \frac{i}{2} \sin \Omega, \\
        s &= \tan \frac{i}{2} \cos \Omega,\\
        \tr{L} &= \Omega + \omega + \nu.
    \end{aligned}
\end{equation}
\tr{The vector $\bm{r}$ describes the shape and orientation of the orbit, while $\tr{L}$ specifies the position along the orbit. The full state vector is given by $\bm{x}=\left[\bm{r}^{\tr{\textrm{T}}}, \tr{L}, m \right]^{\tr{\textrm{T}}}$} and $\bm{u}=\left[u^{\rm R},u^{\rm T},u^{\rm N}\right]$ is the thrust vector in the \gls*{RTN} reference frame. 
\tr{Note that the mean longitude $\tr{L}$ is excluded from the terminal cost and constraints, meaning the final position along the orbit is not considered significant in the optimisation.}
The stage cost is the one from \Eq{eq:fuel_optimal_ctg}, and the path constraints are from \Eq{eq:energy_optimal_problem}.
The normalisation units are reported in \Tab{tab:EARTH_units}, the spacecraft parameters are similar as those of \Tab{tab:SUN_units}, and the tolerances are the same.
\begin{table}[]
    \centering
    \caption{Earth-centered normalisation units and parameters.}
    \begin{tabular}{ l c c}
        \toprule
      \textbf{\tg{Parameter}} & \textbf{Symbol} & \textbf{Value}  \\
      \midrule
        Mass parameter [\SI{}{\kilo\meter\cubed}/\SI{}{\second\squared}] & $ \mu $ & \SI{398600}{} \\
        Time [\SI{}{\second}] & $ \mathrm{TU} $ & \SI{86400}{} \\
        Length [\SI{}{\kilo\meter}] & $ \mathrm{LU} $ & \SI{42241}{}  \\ 
        Velocity [\SI{}{\kilo\meter}/\SI{}{\second}] & $ \mathrm{VU} $ & \SI{0.48890}{} \\
        \bottomrule
    \end{tabular}
    \label{tab:EARTH_units}
\end{table}
Three test cases of the Earth-centered two-body problem will be handled:
\tp{\begin{enumerate}
    \item A transfer from a \gls*{LEO} to another \gls*{LEO} inspired by \cite{DiCarloVasile_2021_ASfLTOT}
    \item A transfer from a \gls*{MEO} to another \gls*{MEO} with a \SI{-35}{\deg} change in $\Omega$
    \item A \gls*{GTO} to \gls*{GEO} transfer from \citet{YangEtAl_2016_OLTSTULBG}.
\end{enumerate}}
\begin{table}[]
    \centering
    \caption{Earth-centered two-body problem transfers data.}
    \begin{tabular}{l c c c c c c c c c} 
        \toprule
        \textbf{Test case} & \textbf{\gls*{ToF}} [\SI{}{\day}] & $N$ & \textbf{Type} & \textbf{$a$ [\SI{}{\kilo\meter}]} & \textbf{$e$ [-]} & \textbf{$i$ [\SI{}{\deg}]} & \textbf{$\Omega$ [\SI{}{\deg}]} & \textbf{$\omega$ [\SI{}{\deg}]} & \textbf{$\nu$ [\SI{}{\deg}]}  \\ [0.5ex] 
        \midrule
        \multirow{2}{*}{\parbox{2cm}{\gls*{LEO} to \gls*{LEO}}} & \multirow{2}{*}{\num{35}} & \multirow{2}{*}{\num{1000}} & Initial & \num{6778.0} & \num{0} & \num{51} & \num{145} & – & \num{0} \\
        & & & Target & \num{7178.0} & \num{0} & \num{56} & \num{145} & – & –  \\
        \midrule[0.5pt]
        \multirow{2}{*}{\parbox{2cm}{\gls*{MEO} to \gls*{MEO}}} & \multirow{2}{*}{\num{55}} & \multirow{2}{*}{\num{1000}} & Initial & \num{34378} & \num{0} & \num{60} & \num{180} & – & \num{0} \\
        & & & Target & \num{34378} & \num{0} & \num{60} & \num{155} & – & –  \\
        \midrule[0.5pt]
        \multirow{2}{*}{\parbox{2cm}{\gls*{GTO} to \gls*{GEO}}} & \multirow{2}{*}{\num{90}} & \multirow{2}{*}{\num{1200}} & Initial & \tp{\num{24505.9}} & \tp{\num{0.725}} & \tp{\num{7}} & – & \num{0} & \num{0}\\
        & & & Target & \num{42165} & \num{0} & \num{0} & – & – & – \\
        \bottomrule
    \end{tabular}
    \label{tab:ic_E_tbp}
\end{table}

\begin{figure}[H]
    \centering 
      \includegraphics[width=0.66\linewidth]{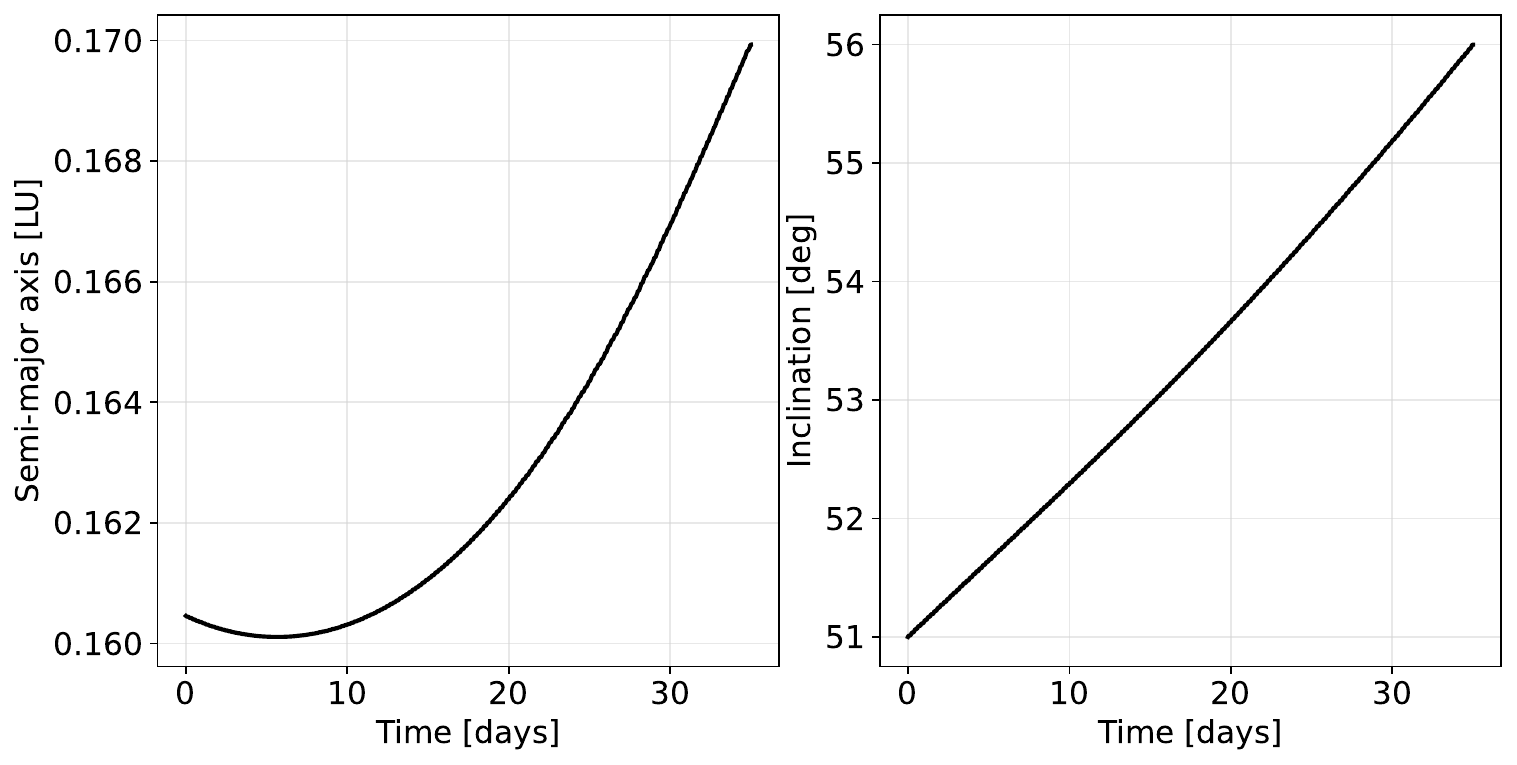}
    \caption{Solution to the \gls*{LEO} to \gls*{LEO} transfer.}
    \label{fig:leo_to_leo}
\end{figure}

\begin{figure}[H]
    \centering 
      \includegraphics[width=\linewidth]{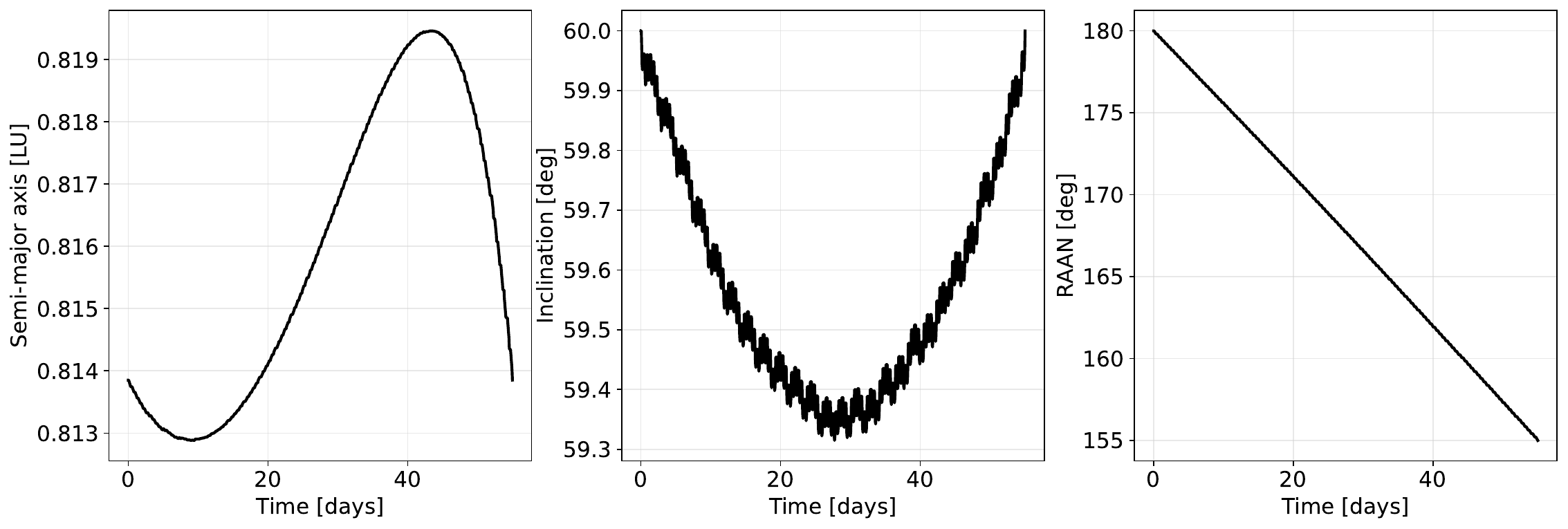}
  
    \caption{Solution to the \gls*{MEO} to \gls*{MEO} transfer.}
    \label{fig:meo_to_meo}
\end{figure}

\begin{figure}[H]
    \centering 
      \includegraphics[width=\linewidth]{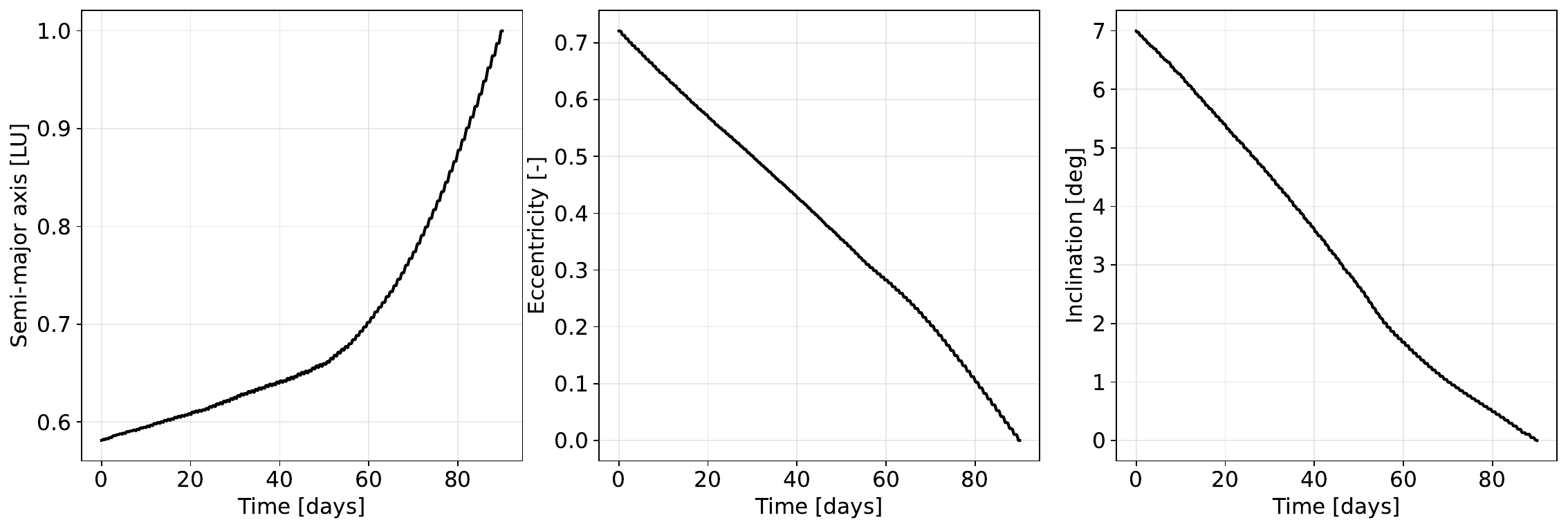}
    \caption{Solution to the \gls*{GTO} to \gls*{GEO} transfer.}
    \label{fig:gto_to_geo}
\end{figure}

\Fig{fig:leo_to_leo} shows \tg{the evolution of the Keplerians}. 
Similarly for \Fig{fig:meo_to_meo} for the \gls*{MEO} to \gls*{MEO} transfer.
\tp{Finally, \Fig{fig:gto_to_geo} presents the \gls*{GTO} to \gls*{GEO} transfer.}
The solutions are also bang-bang controls, the performance metrics are given in \Tab{tab:DADDy_res_E_tbp}, where $n$ DDP stands for the number of DDP iterations performed.
\begin{table}[]
    \centering
    \caption{Performance metrics for fuel-optimal Earth-centred two-body problem transfers.}
    \begin{tabular}{l l c c c c c c} 
        \toprule
        \textbf{Transfer} & \textbf{Data} & \textbf{ iLQR} & \textbf{ \gls*{DDP}} & \textbf{Q} & \textbf{ iLQRDyn} & \textbf{ \gls*{DDP}Dyn} &\textbf{QDyn} \\ [0.5ex] 
        \midrule
    
        \multirow{3}{*}{{\parbox{1cm}{\gls*{LEO} to \gls*{LEO}}}} & $J$ [\SI{}{\kilo\gram}] & \num{40.61} & \num{38.8} & \num{39.77} & \num{40.41} & \num{38.91} & \multirow{3}{*}{\Gls*{DNC}} \\
        & $n$ \gls*{DDP} & \num{765} & \num{696} & \num{439}  & \num{955} & \num{918} & \\
        & Run time [\SI{}{\minute}] & \num{59.4} & \num{56.4} & \num{34.8}  & \num{28.2}& \num{26.1}  & \\
        \midrule[0.5pt]
        \multirow{3}{*}{{\parbox{1cm}{\gls*{MEO} to \gls*{MEO}}}} & $J$ [\SI{}{\kilo\gram}] & \num{74.32} & \num{74.76} & \multirow{3}{*}{\Gls*{DNC}} & \num{74.30} & \num{74.31} & \multirow{3}{*}{\Gls*{DNC}} \\
        & $n$ \gls*{DDP} & \num{1314} & \num{755} &  & \num{2200} & \num{1302}  & \\
        & Run time [\SI{}{\minute}] & \num{59.3}  & \num{28.5} &  & \num{35.3}  & \num{20.3}  & \\
        \midrule[0.5pt]
        \multirow{3}{*}{{\parbox{1cm}{\gls*{GTO} to \gls*{GEO}}}} & $J$ [\SI{}{\kilo\gram}] &  \num{88.5} & \multicolumn{2}{c}{\multirow{3}{*}{\Gls*{DNC}}} & \num{89.18} &     \multicolumn{2}{c}{\multirow{3}{*}{\Gls*{DNC}}}   \\
        & $n$ \gls*{DDP}  & \num{2298} &  &  & \num{3150} & &  \\
        & Run time [\SI{}{\hour}] & \num{3.70} &  & & \num{1.66}  &   &  \\
        \bottomrule
    \end{tabular}
    \label{tab:DADDy_res_E_tbp}
\end{table}
\tp{When they converge, \Gls*{DDP} and Q methods exhibit significantly faster convergence than the \gls*{iLQR} methods. For these multi-revolution transfers, variations of a few percent on the fuel consumption can be observed, especially for the LEO-to-LEO transfer.
Similarly to test cases presented earlier in this work, methods with polynomial dynamics approximation always converge faster than methods than their counterpart that recompute the complete dynamics at each stage.
Note that the run times are rather long, ranging from \SI{26}{\minute} to \SI{59}{\minute} for transfer~\num{1}, from \SI{20}{\minute} to \SI{59}{\minute} for transfer~\num{2}, and from \SI{1.7}{\hour} to \SI{3.7}{\hour} for transfer~\num{3}.
These performances can be explained by the fact that these problems require numerous switches and that the dynamics are implemented without any form of regularisation, such as the Sundman transform \cite{AzizEtAl_2018_LTMRTOVDDPaaST}, to improve convergence. However, these test cases provide a valuable stress test for the \gls*{DADDy} solver and demonstrate its ability to handle complex transfers.}
Averaged analytical methods are a fast and adapted alternative to solve similar optimisation problems \cite{DiCarloVasile_2021_ASfLTOT,DiCarloEtAl_2021_EAFftPKMuLTAaOP}.

\section{Conclusions} \label{sec:conclusions}
\tp{In this work we propose an accelerated approach for constrained spacecraft trajectory optimisation. Building on existing methods, we leverage high‑order Taylor expansions for both automatic differentiation and non-linear‑dynamics approximation. The resulting publicly available \gls*{DADDy} solver integrates a \gls*{DDP}/\gls*{iLQR} routine to generate an optimal, nearly feasible solution without requiring a good initial guess. An enhanced Newton solver then enforces full feasibility, dramatically speeding up the overall constrained‑optimisation process. Compared with the current state‑of‑the‑art, this algorithm delivers substantial runtime reductions across multiple benchmark problems, underscoring its novel contributions and performance gains.}


\tr{The use of high-order Taylor polynomials provides two key advantages: first, it enables automatic differentiation, allowing users to optimise arbitrary dynamical systems with general cost functions and constraints without the need to manually derive gradients and Hessians; second, it enables fast polynomial-based approximations of non-linear dynamics, significantly reducing the computational burden of repeated function evaluations during optimisation \tp{and solution polishing}.}

\tr{While \gls*{DDP} is traditionally more computationally demanding than \gls*{iLQR} due to the need for second-order derivatives of the dynamics, the \gls*{DA} framework mitigates this overhead. By computing derivatives and generating polynomial approximations simultaneously, the runtime per iteration of \gls*{DDP}Dyn (i.e., \gls*{DDP} with dynamics approximation) is brought closer to that of \gls*{iLQR}Dyn (\gls*{iLQR} with approximation of the dynamics).
Results show that \num{79.4}\% of dynamics evaluations are handled via high-order Taylor approximations, substantially accelerating the overall process. This framework also enabled the implementation of a "Q" method that directly evaluates the cost-to-go function and its derivatives.}

\tr{Experimental results on various test cases show that the \gls*{iLQR}Dyn method is the most stable among the tested optimisation methods, achieving results comparable to those in the literature while running \tp{\num{41}\% to \num{88}\%} faster than the standard \gls*{iLQR} method, with no observed drawbacks.
The \gls*{DDP} and \gls*{DDP}Dyn methods outperform \gls*{iLQR} and \gls*{iLQR}Dyn in terms of runtime when they converge, confirming the faster convergence rate of \gls*{DDP}.
However, they do not always converge and may require additional regularisation to match the robustness of \gls*{iLQR}-based methods.
The Q and QDyn methods, while mathematically equivalent to \gls*{DDP} and \gls*{DDP}Dyn, exhibit higher run times and greater instability.}

\tr{Finally, the set of \gls*{DADDy} methods performs well across a wide range of trajectory optimisation problems, achieving low run times and satisfactory constraint satisfaction.}
\tp{The algorithm can also optimise many‑revolution transfers with numerous stages, such as Earth‑centred low‑thrust trajectories.}

\section*{Funding sources}
This work was funded by SaCLaB (grant number 2022-CIF-R-1), a research group of ISAE-SUPAERO. 
\section*{Acknowledgments}
The authors acknowledge the use of an AI tool (Claude Haiku 4.5) to assist in the reformulation of certain sentences.

\section*{References}
\bibliographystyle{unsrtnat}
\bibliography{references}

\end{document}